\documentclass[a4paper,11pt]{article}
\usepackage[utf8]{inputenc}
\usepackage[T1]{fontenc}
\usepackage{amsmath,amsfonts,amssymb}
\usepackage{fullpage}
\usepackage[english]{babel}
\usepackage{graphicx}
\usepackage{titlesec}
\usepackage[runin]{abstract}
\usepackage{enumitem}
\usepackage{latexsym}
\usepackage{color}
\usepackage{amsmath}
\usepackage{mathtools,amscd}
\usepackage{ae, aeguill}
\usepackage{latexsym}
\usepackage[all]{xy}
\renewcommand{\a }{\alpha }
\renewcommand{\b }{\beta }
\renewcommand{\d}{\delta }

\newcommand{\D }{\Delta }

\newcommand{\e }{\varepsilon }
\newcommand{\g }{\gamma}

 \newcommand{\foral }{\forall\, }
\newcommand{\G }{\Gamma }
\renewcommand{\l }{\lambda }

\newcommand{\m }{\mu }
\newcommand{\na }{\nabla }

\newcommand{\rh }{\rho }
\newcommand{\s }{\sigma }
\newcommand{\Sig }{\Sigma}

\renewcommand{\th }{\theta }
\renewcommand{\o }{\omega }
\renewcommand{\O }{\Omega }

\newcommand{\ov}{\overline}
\newcommand{\wtilde }{\widetilde}

\newcommand{\pa}{\partial}
\newcommand{\beq}{\begin{eqnarray}}
\newcommand{\eeq}{\end{eqnarray}}
\newcommand{\beqq}{\begin{eqnarray*}}
\newcommand{\eeqq}{\end{eqnarray*}}
\newcommand{\be}{\begin{equation}}
\newcommand{\ee}{\end{equation}}
\newcommand{\ba}{\begin{align*}}
\newcommand{\ea}{\end{align*}}
\newenvironment{pf}{\noindent{\bf Proof.}\enspace}{
\hfill$\Box$\medskip}

\newcommand{\Real}{\mathbb{R}}
\newcommand{\norm}[1]{\Vert {#1} \Vert} 

\newcommand{\diff}{\mathrm{d}}%
\newcommand{\soexp}[1]{{#1}_{\m}^{\ast}}
\renewcommand{\norm}[1]{\left\Vert#1\right\Vert}

\newcommand{\abs}[1]{\left\vert#1\right\vert}
\newcommand{\set}[1]{\left\{#1\right\}}
\newcommand{\seq}[1]{\left<#1\right>}

\usepackage{tikz}
\usetikzlibrary{tikzmark,calc}

\newcommand{\tikzarc}[1]{%
  \tikzmarknode{a}{#1}%
  \begin{tikzpicture}[overlay, remember picture]
    \draw ([yshift=1pt]a.north west) to[bend left=20] ([yshift=1pt]a.north east);
    \node at ($(a.north)!0.01!(a.north east) + (0, 8pt)$) {°};
  \end{tikzpicture}%
}

\newcommand{\Z}{\mathbb{Z}}

\newcommand{\Natl}{\mathbb{N}}

\newtheorem{theorem}{Theorem}[section]

\newtheorem{lemma}{Lemma}[section]
\newtheorem{proposition}{Proposition}[section]

\numberwithin{equation}{section}
\let \n = \noindent

\title{
\begin{center}
\vspace{0.5in}{\bf\Large Nonlocal Choquard equations\\ involving critical Hardy-Littlewood-Sobolev\\exponent: the effect of the topology of the domain}
\end{center}}
\date{}
\begin{document}
\maketitle
\begin{center}
\author{{\bf\large Mohammed Ali Mohammed Alghamdi}} 
{\bf\large}\vspace{1mm}\\
{\it\small Department of Mathematics},\\
{\it\small King Abdulaziz University},
{\it\small  Jeddah 21589, Saudi Arabi}\\
{\it\small mamalgamdi@kau.edu.sa}\vspace{1mm}\\
{\bf\large Hichem Chtioui } 
{\bf\large}\vspace{1mm}\\
{\it\small Sfax University},\\ {\it\small
Faculty of Sciences of Sfax},\\
{\it\small  3018 Sfax, Tunisia.}\\
{\it\small Hichem.Chtioui@fss.rnu.tn }\vspace{1mm}\\

\end{center}
\begin{abstract}. We apply a topological method to prove existence of positive solutions for the nonlinear Choquard equation with upper critical exponent in the sense of Hardy-Littlewood-Sobolev inequality on bounded domains having nontrivial homology group. \\

\n Keywords: Choquard equation, Critical nonlinearity, Variational calculus, Singular homology, Topological method.
\end{abstract}
\section{Introduction}\label{sec1}
During the last decades, the Choquard equation
\beq\label{eq1.1}
  -\Delta\,u+V(x)u &=& \left(\int_{\Real^n}\frac{\abs{u(y)}^q}{\abs{x-y}^\mu}\diff y\right)\abs{u}^{q-2}u \quad\text{in }\Real^n,
\eeq
where $n\geq 3$,\;$0<\mu<n$ and $\frac{2n-\mu}{n}\leq q\leq\frac{2n-\mu}{n-2}$, has attracted the attention and the interest of a lot of researchers, since it appears in the modeling of various physical phenomena, such as quantum mechanics of a polaron at rest \cite{29 BN}, Selfgravitatinal collapse of a quantum mechanical wave function \cite{30 BN}, modeling an electron trapped in its own hole as an approximation to Hartree-Fock theory of one component plasma \cite{21 BN} and so on.

\n Following standard variational theory, the solutions of the Choquard equation are the critical points of the Euler-Lagrange functional
\beqq
I(u) &=& \frac{1}{2}\int_{\Real^n}\abs{\na u}^2\diff x-\frac{1}{2q}\int_{\Real^n}\int_{\Real^n}\frac{\abs{u(y)}^q\abs{u(x)}^q}{\abs{x-y}^\m}\diff x\diff y,
\eeqq
on $H^1(\Real^n)$. Recall that by the celebrate Hardy-Littlewood-Sobolev inequality see \cite{15 DCDS} and \cite{17 DCDS}, the second term of $I(u)$ is well defined on $H^1(\Real^n)$, if and only if
\beqq
2_\m:=\frac{2n-\m}{n} &\leq q\leq &\frac{2n-\m}{n-2}:=2_{\m}^{\ast}.
\eeqq
In this case, the functional $I$ is continuously differentiable on $H^1(\Real^n)$. The constants $2_\m$ and $\soexp{2}$ are termed, respectively, as the lower critical exponent and the upper critical exponent in the sense of Hardy-Littlewood-Sobolev. It is remarkable that the Choquard equation is a nonlinear problem with superlinear nonlinearity, since $2\m>1$, $\foral\m\in(0,n)$, and it is a nonlocal problem. The nonlocal term which appears in the nonlinearity makes the problem particularly difficult. For $n=3$, $q=2$, $\m=1$ and $V$ is a positive constant, Lieb \cite{21 BN} obtained the existence of a unique minimizing positive solution (up to translations). Under some restrictions on $V(x)$, Lions \cite{19 DCDS} proved the existence of a sequence of radially solutions. For $n\geq 3$, $\m\in(0,n)$, $V=1$, and suitable range of $q$, Moros-Van Schaftingen, (
see \cite{19 2023} and \cite{23 2023}), proved the existence of a nontrivial solution of \eqref{eq1.1}, if and only if $2\m<q<\soexp{2}$. Their result has obtained by using a Pohozaev type identity and minimizing arguments. Ma-Zhao \cite{25 DCDS} studied the positive and regularity of the minimizing solutions. For $q=2_\m$, the lower critical exponent and a nonconstant potential $V$, Moros-Van Schaftingen \cite{47 Aut} and Cassani-Schaftingen-Zhang \cite{16 Aut} studied the existence problem of positive solutions of equation \eqref{eq1.1}. they provided sufficient and necessary conditions on the potential $V(x)$ to obtain minimizing solutions. For more related results on Choquard type equations, we refer to \cite{P2}, \cite{P1},
\cite{P3} and references therein.

\n For the upper critical exponent $q=\soexp{2}$ and $V=0$, the Choquard equation reduces to 
\beq\label{eq1.2}
-\D u &=& \left(\int_{\Real^n}\frac{\abs{u(y)}^{\soexp{2}}}{\abs{x-y}^\m}\diff y\right)\abs{u}^{2_\m^{\ast}-2}u\quad\text{in }\Real^n.
\eeq
Using moving plane techniques, Du-Yang \cite{17 Pd} and Guo-Hu-Peng-Shuai \cite{23 Pd} proved a uniqueness classification of the positive solutions of \eqref{eq1.2}. Precisely, they proved that the positive solutions of \eqref{eq1.2} have to be of the form
\beqq
\wtilde\d_{(a,\l)}(x) &=& (n(n-2))^{\frac{n-2}{4}}C_{n,\m}^{\frac{2-n}{2(n-\m+2)}}S^{\frac{(n-\m)(2-n)}{4(n-\m+2)}}\d_{(a,\l)}(x),\;\;x\in\Real^n,
\eeqq
where
\beqq
C_{n,\m} &=& \pi^{\frac{\m}{2}}\frac{\G(\frac{n-\m}{2})}{\G(n-\frac{\m}{2})}(\frac{\G(n)}{\G(\frac{n}{2})})^{\frac{n-\m}{n}}, \quad\G(s)=\int_0^{\infty}\frac{e^{-t}}{t^{1-s}}\diff t, \,s>0,
\eeqq
is the Gamma function, $S$ is the best constant of Sobolev and for $a\in\Real^n$ and $\l>0$,
\beqq
\d_{(a,\l)}(x) &=& \frac{\l^{\frac{n-2}{2}}}{(1+\l^2\abs{x-a}^2)^{\frac{n-2}{2}}}\cdot
\eeqq
In this paper we restrict our attention to problem \eqref{eq1.2} on bounded domains $\O$ of $\Real^n$, $n\geq 3$. we are interested to study the existence problem of positive solutions of the following nonlocal Choquard equation with upper critical Hardy-Littlewood-Sobolev exponent
\beq\label{eq1.3}
-\D u &=& \left(\int_{\O}\frac{\abs{u(y)}^{\soexp{2}}}{\abs{x-y}^\m}\diff y\right)\abs{u}^{2_\m^{\ast}-2}u\quad\text{in }\O,\;\;u=0\text{  on  }\partial\O,
\eeq
where $\m\in(0,n)$ and $2_\m^{\ast}=\frac{2n-\m}{n-2}$.

\n Since $\soexp{2}$ is critical, a lack of compactness occurs in the sense that the associated Euler Lagrange functional does not satisfy the Palais-Smale condition. This makes the problem of finding positive solutions more challenging. As the celebrate critical elliptic local equation
\beq\label{eq1.4}
-\D u &=& \abs{u}^{2^{\ast}-2}u\quad\text{in }\O,\quad u=0\text{  on  }\pa\O,
\eeq
\n where $2^{\ast}=\frac{2n}{n-2}$ is the critical Sobolev exponent, a concentration compactness principle \cite{Lions}, \cite{Struwe} and Bahri-Coron representations \cite{BC 1991}, \cite{BC 1989} play an important role in solving nonlinear problems of a non compact nature.

\n In \cite{31 Aut}, Goel-Radulescu-Sreenadh considered problem \eqref{eq1.3} where $\O$ is a bounded annular-type domain of $\Real^n,\;\;n\geq 3$. They proved a global compactness Lemma and used a variational method to show existence of positive solutions when the inner hole of $\O$ is sufficiently small. Their result is motivated by the celebrate paper of Coron \cite{C} for the local equation \eqref{eq1.4}. In \cite{20 Pd}, Gao-Yang studied the Hartree type Brezis-Nirenberg problem; that is a perturbation of equation \eqref{eq1.3} by adding a linear term $"\l u,\;\l\in\Real"$ to the Choquard nonlinearity of \eqref{eq1.3}. They proved existence results for a suitable range of $\l$ and $a$ non-existence result in the case of $\l\leq 0$ and $\O$ is star shaped domain with respect to the origin. In the papers \cite{20 Pd}, \cite{19 Green}, \cite{11 2023} and \cite{SYZ Pdelta}, the method of Brezis-Nirenberg developed for the local elliptic equation, see \cite{BN}, has been adopt to study nonlocal Choquard problems with upper critical nonlinearity.

\n Pushing more the resemblance of the local equation \eqref{eq1.4} and the nonlocal Choquard equation \eqref{eq1.3}, we are led to investigate the effect of the domain's topology on the existence results of problem \eqref{eq1.3}. In the pioneering paper \cite{BC 1989}, Bahri-Coron introduced a powerful topological argument which leads to prove the existence of positive solutions of equation \eqref{eq1.4} when $\O$ has a nontrivial homology group. In this paper we prove that the topological approach of Bahri-Coron can be successfully adopt to study the existence problem of positive solutions of the Choquard equation \eqref{eq1.3}. Because of its nonlocal nature, additional difficulties arise in extending such a topological method to the present framework. This requires changes in the construction and novelties in the proof.

\n Let $k$ be a non-negative integer. We denote $H_k(\O)$ the singular homology group of $\O$ of dimension $k$ with $\Z_2$-coefficients. We shall prove the following result.
\begin{theorem}\label{th1.1}
Let $n\geq 3$ and let $\O$ be a bounded domain of $\Real^n$. If there exists a positive integer $k_0$ such that $H_{k_0}(\O)$ is not null, then problem \eqref{eq1.3} has at least a positive solution.
\end{theorem}
The proof of Theorem \ref{th1.1} will be the subject of Section \ref{sec3} of this paper. We summarize the idea of the proof in the following three steps.

\n For a pair of topological spaces $(X,Y)$, $Y\subset X$, we denote $H_\ell(X,Y)$, $\ell\in\Natl$, the relative homology of $(X,Y)$.
\begin{description}[leftmargin=*,parsep=0.cm,itemsep=0.cm,topsep=0cm,resume]
  \item[$\text{Step}\;1.$] Under the assumption that $H_{k_0}(\O)$ is not trivial, we prove that there exist two sequences of topological pairs $(B_p,B_{p-1})$ and $(E_p,F_{p})$, $p\in\Natl\backslash\{0\}$, a sequence of non-zero homology classes $(\o_p)$ such that

      $\o_p\in H_\ell(B_p,B_{p-1})$, for any $p$ and a sequence of relative homology homomorphisms
      \[(\phi_p)_\ell\colon H_\ell(B_p,B_{p-1})\to H_\ell(E_p,F_p)\]
      such that the first homomorphism of the sequence satisfies
      \beq\label{eq1.5}
      (\phi_1)_\ell(\o_1)&\neq& 0.
      \eeq
  \item[$\text{Step}\;2.$] If we suppose that problem \eqref{eq1.3} has no positive solution, we prove that

      for any $p\in\Natl\backslash\{0\}$, \beq\label{eq1.6}
      (\phi_p)_\ell(\o_p)\neq 0\,&\Longrightarrow&\,(\phi_{p+1})_\ell(\o_{p+1})\neq 0.
      \eeq
  \item[$\text{Step}\;3.$] We prove the existence of a positive large integer $p_0$ such that
      \beq\label{eq1.7}
      (\phi_{p_0})_\ell &=& 0\,\text{  for any  }\ell\in\Natl.
      \eeq
      This concludes the proof.
\end{description}
In order to achieve the above steps, we state in the next Section the variational structure associated to problem \eqref{eq1.3} including some preparatory results. We first change the Euler-Lagrange functional $I$ and work with a more convenient functional $J$. We study the behavior of sequences failing the Palais-Smale condition and Parametrize the neighborhoods of critical points at infinity \cite{BL}. Finally, we prove asymptotic estimates and expand the Euler-Lagrange functional $J$ in these neighborhoods. Particularly, we determine exactly the levels of $J$ where the Palais-Smale condition does not hold.

\n In the following we denote, $\seq{u,v}=\int_{\O}\na u\na v\diff x$ the inner product of $H_0^1(\O)$,

\n $\abs{u}_{1,\O}=\left(\int_\O\abs{\na u}^2\diff x\right)^{\frac{1}{2}}$ the associated norm and the $\|u\|_{L^t}=\left(\int_\O\abs{u}^t\diff x\right)^{\frac{1}{t}}$  norm of $L^t(\O)$.
\section{Variational results and asymptotic expansions}\label{sec2}
\n We start this Section by recalling the Hardy-Littlewood-Sobolev inequality which play a fundamental role in defining the variational approach to critical problems with nonlocal Choquard nonlinearities.
\begin{lemma}\label{lem2.1}
\cite{15 DCDS}, \cite{17 DCDS}. For $t_1,\;t_2>1$, such that $\frac{1}{t_1}+\frac{1}{t_2}=\frac{2n-\m}{n}$, we have
\beq\label{eq2.1}
\int_{\Real^n}\int_{\Real^n}\frac{f(y)g(x)}{\abs{x-y}^\m}\diff x\diff y &\leq& C\norm{f}_{L^{t_1}}\norm{g}_{L^{t_2}},\;\forall(f,g)\in L^{t_1}(\Real^n)\times L^{t_2}(\Real^n).
\eeq
Here $C$ is a positive constant independent of $f$ and $g$.

\n Particularly, if $t_1=t_2=\frac{2n}{2n-\m}$, then $C$ equals to $C_{n,\m}$; the positive constant introduced in the above Section. In this case, there is equality in \eqref{eq2.1}, if and only if $f=g$ (up to a multiplicative constant) and
\beqq
g(x) &=& A(\l^2+\abs{x-a}^2)^{\frac{\m-2n}{2}},
\eeqq
where $A\in\mathbb{C}$, $\l\in\Real\backslash\{0\}$ and $a\in\Real^n$.
\end{lemma}
Let us observe that under the Hardy-Littlewood-Sobolev inequality, the integral
\beq\label{eq2.2}
\int_{\Real^n}\int_{\Real^n}\frac{\abs{u(y)}^q\abs{u(x)}^q}{\abs{x-y}^\m}\diff x\diff y &&
\eeq
is well defined, if $\abs{u}^q\in L^t(\Real^n)$ such that $\frac{2}{t}=\frac{2n-\m}{n}$. Thus, by Sobolev embedding, integral \eqref{eq2.2} is defined for $u\in H^1(\Real^n)$, only if $\frac{2n-\m}{n}\leq q\leq\frac{2n-\m}{n-2}$.

\n In the special case of the upper critical Hardy-Littlewood-Sobolev exponent $2^\ast_\m=\frac{2n-\m}{n-2}$ and for an open set $\O$ of $\Real^n$, $n\geq 3$, it is proved in \cite{20 Pd} that
\beqq
\norm{u}_{HL} &=& \left(\int_{\O}\int_{\O}\frac{\abs{u(y)}^{\soexp{2}}\abs{u(x)}^{\soexp{2}}}{\abs{x-y}^\m}\diff x\diff y\right)^{\frac{1}{2\,\soexp{2}}},\;\;u\in L^{2^\ast}(\O),
\eeqq
is a norm on $L^{2^\ast}(\O)$, where $2^\ast=\frac{2n}{n-2}$. In this case the best constant $S_{HL}(\O)$ defined by
\beq\label{eq2.3}
S_{HL}(\O) &=& \inf_{u\in H_0^1(\O)\setminus\{0\}}\frac{\abs{u}_{1\O}}{\norm{u}_{HL}},
\eeq
is independent of $\O$ and it is never achieved except $\O=\Real^n$, see \cite{20 Pd}. In addition, the unique minimizers of $S_{HL}(\Real^n)$ are unique and of the form $u=c\d_{(a,\l)}$, where $c$ is a positive constant and the parameters $(a,\l)\in\Real^n\times\Real^+$. Here
\beq\label{eq2.4}
\d_{(a,\l)} &=& \frac{\l^{\frac{n-2}{2}}}{(1+\l^2\abs{x-a}^2)^{\frac{n-2}{2}}}\;,\;\;x\in\Real^n.
\eeq
As a consequence of the Hardy-Littlewood-Sobolev inequality, the functional $I$ defined by
\beqq
I(u) &=& \frac{1}{2}\abs{u}_{1\O}^2-\frac{1}{2\,\soexp{2}}\norm{u}_{HL}^{2\,\soexp{2}},
\eeqq
is well defined  on $H_0^1(\O)$ and it is of class $C^1$. It is straightforward to see that for any $u\in H_0^1(\O)$, we have
\beqq
I'(u) &=& -\D u-\left(\int_\O\frac{\abs{u(y)}^{\soexp{2}}}{\abs{x-y}^\m}\diff y\right)\abs{u}^{\soexp{2}-2}u.
\eeqq
Therefore the solutions of problem \eqref{eq1.3} are the critical points of $I$. Since $\soexp{2}$ is critical in the sense of Hardy-Littlewood-Sobolev inequality, the functional $I$ does not satisfy the Palais-Smale condition. A description of sequences failing the Palais-Smale condition has been established in \cite{31 Aut}. As a product, existence of a positive critical point of $I$ is obtained under the assumption that $\O$ is a bounded annular-type domain with inner hole small enough. Let
\beqq
\Sig &=& \set{u\in H_0^1(\O),\;\;\abs{u}_{1\O}=1}.
\eeqq
On $\Sig$, we define
\beqq
J_1(u) &=& \sup_{\l\geq 0}I(\l u).
\eeqq
\begin{lemma}\label{lem2.2}
For any $u\in\Sig$, there exists a unique $\l=\l(u)>0$ such that $J_1(u)=I(\l(u)u)$.
\end{lemma}
\begin{pf}
Let $u\in\Sig$ and let $f(\l)=I(\l u)$, $\l\geq 0$. We have
\beqq
f'(\l) &=& \l-\l^{2\,\soexp{2}-1}\norm{u}_{HL}^{2\,\soexp{2}}.
\eeqq
Therefore, $f$ is an increasing function near $0$, $f(0)=0$ and $f(\l)\to -\infty$, as $\l\to +\infty$. Then there exists $\l=\l(u)>0$ such that $f(\l(u))=\underset{\l\geq0}{\max}f(\l)$. Since $f'(\l(u))=0$, then $\l(u)$ is unique and equals to $\norm{u}_{HL}^{\frac{-\soexp{2}}{\soexp{2}-1}}$.
\end{pf}

\n Let us observe that there is a correspondence between the critical points of $J_1$ and the critical points of $I$. Indeed,
\begin{lemma}\label{lem2.3}
Let $u\in\Sig$. $u$ is a critical point of $J_1$ if and only if $\l(u)u$ is a critical point of $I$.
\end{lemma}
\begin{pf}
For any $u\in\Sig$, we have by Lemma \ref{lem2.2}.
\beqq
J'_1(u) &=& \l(u)I'(\l(u)u)+I'(\l(u)u)(u).\l'(u).
\eeqq
Using the fact that $f'(\l(u))=0$, where $f(\l)$ is defined in the proof of Lemma \ref{lem2.2}, we get $I'(\l(u)u)(u)=0$ and hence $J'_1(u)=I'(\l(u)u)$, up to a positive multiplicative constant. The result follows.
\end{pf}

\n Note that by the result of Lemma \ref{lem2.2}, $J_1(u)$, $u\in\Sig$, can be expressed by
\beqq
J_1(u) &=& (\frac{1}{2}-\frac{1}{2\,\soexp{2}})\frac{1}{\left(\int_\O\int_\O\frac{\abs{u(y)}^{\soexp{2}}\abs{u(x)}^{\soexp{2}}}{\abs{x-y}^\m}\diff x\diff y\right)^{\frac{1}{\soexp{2}-1}}},
\eeqq
since $\l(u)=\norm{u}^{-\frac{\soexp{2}}{\soexp{2}-1}}$. Moreover, by Lemma \ref{lem2.3}, $J(u):=(\frac{1}{2}-\frac{1}{2\,\soexp{2}})^{-1}J_1(u)^{\soexp{2}-1}$, can be also considered as an Euler-Lagrange functional associated to problem \eqref{eq1.3}, in the sense that the solutions of \eqref{eq1.3} correspond to critical points of $J$. Then, to prove our result, it is more convenient for us to work with,
\beqq
J(u) &=& \frac{1}{\int_\O\int_\O\frac{\abs{u(y)}^{\soexp{2}}\abs{u(x)}^{\soexp{2}}}{\abs{x-y}^\m}\diff x\diff y},\;\;u\in\Sig.
\eeqq
By Hardy-Littlewood-Sobolev inequality and \eqref{eq2.3}, the functional $J$ is lower bounded on $\Sig$ and we have
\beqq
\inf_{u\in\Sig}J(u) &=& S_{HL}(\O)^{2\,\soexp{2}}.
\eeqq
Denote $\wtilde{S}_{HL}=S_{HL}(\O)^{2\,\soexp{2}}$. Using the result of \cite{20 Pd}, $\wtilde{S}_{HL}$ is independent of $\O$ and it is never achieved except $\O=\Real^n$. In this case, the minimizers functions are unique and of the form $c\d_{(a,\l)}$, where $\d_{(a,\l)}$ is defined in \eqref{eq2.4} and $c>0$.

\n Let $a\in\Real^n$ and $\l>0$, we define
\beqq
U_{(a,\l)}(x) &=& \g_0\d_{(a,\l)}(x),\;\;x\in\Real^n,
\eeqq
where
\beq\label{eq2.5}
\g_0=\left(\int_{\Real^n}\abs{\na\d_{(a,\l)}}^2\diff x\right)^{-\frac{1}{2}} &=& \left(n(n-2)\int_{\Real^n}\frac{\diff x}{(1+\abs{x}^2)^n}\right)^{-\frac{1}{2}}.
\eeq
\begin{lemma}\label{lem2.4}
For $a\in\Real^n$ and $\l>0$, $U_{(a,\l)}$ satisfies
\beqq
-\D U_{(a,\l)} &=& \wtilde{S}_{HL}\left(\int_{\Real^n}\frac{U_{(a,\l)}^{\soexp{2}}(y)}{\abs{x-y}^\m}\diff y\right)U_{(a,\l)}^{\soexp{2}-1}\;\;\;\text{in  }\Real^n.
\eeqq
\end{lemma}
\begin{pf}
From \cite{17 Pd} and \cite{20 Pd}, we know that $\d_{(a,\l)}$ satisfies
\beqq
-\D\d_{(a,\l)} &=& A_{HL}\left(\int_{\Real^n}\frac{\d_{(a,\l)}^{\soexp{2}}(y)}{\abs{x-y}^\m}\diff y\right)\d_{(a,\l)}^{\soexp{2}-1}\;\;\;\text{in  }\Real^n,
\eeqq
where, $A_{HL}=(n(n-2))^{\frac{n-\m+2}{2}}C_{n,\m}^{-1}S^{\frac{\m-n}{2}}$. Therefore,
\beqq
-\D U_{(a,\l)} &=& \g_0^{2-2\,\soexp{2}} A_{HL}\left(\int_{\Real^n}\frac{U_{(a,\l)}^{\soexp{2}}(y)}{\abs{x-y}^\m}\diff y\right)U_{(a,\l)}^{\soexp{2}-1}.
\eeqq
Multiplying by $U_{(a,\l)}$ and integrate, we get
\beqq
1 &=& \g_0^{2-2\,\soexp{2}} A_{HL}\int_{\Real^n}\int_{\Real^n}\frac{U_{(a,\l)}(y)^{\soexp{2}}U_{(a,\l)}^{\soexp{2}}(x)}{\abs{x-y}^\m}\diff x\diff y,
\eeqq
and hence
\beqq
\g_0^{2-2\,\soexp{2}} A_{HL} &=& \frac{1}{\int_{\Real^n}\int_{\Real^n}\frac{U_{(a,\l)}(y)^{\soexp{2}}U_{(a,\l)}(x)^{\soexp{2}}}{\abs{x-y}^\m}\diff x\diff y}=\wtilde{S}_{HL}.
\eeqq
\end{pf}

\n For $a\in\O$, we denote
\beqq
G(a,x) &=& \frac{\g_0}{\abs{x-a}^{n-2}}-H(a,x),
\eeqq
where $\g_0$ is defined in \eqref{eq2.5} and
\[\left\{
    \begin{array}{ll}
      \D_x H(a,x) &=0 \hspace{1.8cm}\text{in  }\O\hbox{,} \\
      H(a,x) &=\frac{\g_0}{\abs{x-a}^{n-2}}\hspace{0.5cm}\text{on  }\pa\O \hbox{.}
    \end{array}
  \right.
\]
Our aim in the remainder part of this Section is to establish an asymptotic expansion of the functional $J$ on a suitable finite dimensional set in terms of the Green function $G$ and its regular part $H$.

\n Denote $PU_{(a,\l)}$, $a\in\O$ and $\l>0$, the projection of $U_{(a,\l)}$ on $H_0^1(\O)$. Namely, $PU_{(a,\l)}$ is the unique solution of
\[\left\{
    \begin{array}{ll}
      -\D PU_{(a,\l)} &= -\D U_{(a,\l)} \hspace{1.cm}\text{in  }\O\hbox{,} \\
      PU_{(a,\l)} &=0\hspace{1.5cm}\text{on  }\pa\O \hbox{.}
    \end{array}
  \right.\]
We claim that
\beq\label{eq2.6}
PU_{(a,\l)} &=& U_{(a,\l)}-\frac{H(a,.)}{\l^{\frac{n-2}{2}}}+O\left( \frac{1}{\l^{\frac{n+2}{2}}\diff(a,\pa\O)^n}\right).
\eeq
Indeed, denote $\th_{(a,\l)}=U_{(a,\l)}-PU_{(a,\l)}-\frac{H(a,.)}{\l^{\frac{n-2}{2}}}$. We have
\[\left\{
    \begin{array}{ll}
      \D\th_{(a,\l)} &=0 \hspace{3.8cm}\text{in  }\O\hbox{,} \\
      \th_{(a,\l)}(x) &=U_{(a,\l)}-\frac{\g_0}{\l^{\frac{n-2}{2}}\abs{x-a}^{n-2}}\hspace{0.5cm}\text{on  }\pa\O \hbox{.}
    \end{array}
  \right.
\]
On the boundary $\pa\O$, we have
\beqq
\th_{(a,\l)}(x) &=& \frac{-\g_0}{\l^{\frac{n-2}{2}}\abs{x-a}^{n-2}}\left(1-\biggl(\frac{\l^{2}\abs{x-a}^{2}}{1+\l^{2}\abs{x-a}^{2}}\biggr)^{\frac{n-2}{2}}\right),
\eeqq
and hence
\beqq
\abs{\th_{(a,\l)}(x)} &\leq& \frac{\g_0}{\l^{\frac{n+2}{2}}\abs{x-a}^{n}},\quad\forall x\in\pa\O.
\eeqq
Claim \eqref{eq2.6} follows from the maximum principle Theorem.

\n Let $M$ be a compact set included in $\O$. We then have the following expansion of $J$.
\begin{proposition}\label{prop2.5}
Let $p\geq 1$ and let $a_1,\cdots,a_p\in M$, $\a_1,\cdots,\a_p\geq 0$ such that $\underset{i=1}{\overset{p}{\sum}}\a_i=1$ and $\l>0$. Denote $\diff_a=\underset{i\neq j}{\min}\diff(a_i,a_j)$. If $\l\diff_a$ is large enough, we have
\beqq
J\left(\frac{\sum_{i=1}^{p}\a_i PU_{(a_i,\l)}}{\abs{\sum_{i=1}^{p}\a_i PU_{(a_i,\l)}}_{1\O}}\right) &=& \wtilde{S}_{HL}\frac{(\sum_{i=1}^{p}\a_i^2)^{\soexp{2}}}{\sum_{i=1}^{p}\a_i^{2\,\soexp{2}}}\Biggl[1-\soexp{2}n(n-2)\frac{\g_0 c_1}{\l^{n-2}}\Biggl(\sum_{i=1}^{p}\biggl(\frac{\a_i^2}{\sum_{j=1}^{p}\a_j^{2}}\\
&\hspace{-7.5cm}-& \hspace{-4cm}2\frac{\a_i^{2\, \soexp{2}}}{\sum_{j=1}^{p}\a_j^{2\,\soexp{2}}}\biggr)H(a_i,a_i)
+ \sum_{k\neq i}\biggl(\frac{2\a_i^{2\, \soexp{2}-1}\a_k}{\sum_{j=1}^{p}\a_j^{2\,\soexp{2}}} - \frac{\a_i\a_k}{\sum_{j=1}^{p}\a_j^{2}}\biggr)G(a_i,a_k)\Biggr)\Biggr]+O\biggl(\frac{1}{(\l\diff_a)^{n-1}}\biggr).
\eeqq
Here $\g_0$ is defined in \eqref{eq2.5} and $c_1=\int_{\Real^n}\frac{\diff x}{(1+\abs{x}^2)^{\frac{n+2}{2}}}$.
\end{proposition}
\begin{pf}
Denote $u=\underset{i=1}{\overset{p}{\sum}}\a_i PU_{(a_i,\l)}$.
\beqq
J(\frac{u}{\abs{u}_{1\O}})=\frac{\abs{u}_{1\O}^{2\,\soexp{2}}}{\int_\O\int_\O\frac{u(y)^{\soexp{2}}u(x)^{\soexp{2}}}{\abs{x-y}^\m}\diff x\diff y} &:=& \frac{N}{D}.
\eeqq
We expand the numerator as follows
\beqq
N^{\frac{1}{\soexp{2}}} &=& \sum_{i=1}^{p}\a_i^2\abs{PU_{(a_i,\l)}}_{1\O}^2+\sum_{i\neq j}\a_i\a_j\seq{PU_{(a_i,\l)},PU_{(a_j,\l)}}.
\eeqq
Using estimate \eqref{eq2.6}, for any $i=1,\cdots,p$, we have
\beq\label{eq2.7}
\abs{PU_{(a_i,\l)}}_{1\O}^2 &=& \int_\O-\D U_{(a_i,\l)}PU_{(a_i,\l)}\diff x\nonumber\\
&=& \int_\O-\D U_{(a_i,\l)}\bigg(U_{(a_i,\l)}-\frac{H(a_i,x)}{\l^{\frac{n-2}{2}}}\biggr)\diff x+O\biggl(\frac{1}{\l^{\frac{n+2}{2}}}\int_{\Real^n}U_{(a_i,\l)}^{2^{\ast}-1}\diff x\biggr).
\eeq
Observe that,
\beq\label{eq2.8}
\int_\O-\D U_{(a_i,\l)}U_{(a_i,\l)}\diff x &=&\wtilde{S}_{HL}\int_{\Real^n}\int_{\Real^n}\frac{U_{(a_i,\l)}^{\soexp{2}}(y)U_{(a_i,\l)}^{\soexp{2}}(x)}{\abs{x-y}^\m}\diff x\diff y+O(\frac{1}{\l^n})\nonumber\\
&=& 1+O(\frac{1}{\l^n}).
\eeq
By expanding $H(a_i,x)$ around $a_i$ we have
\beqq
\int_\O-\D U_{(a_i,\l)}\frac{H(a_i,x)}{\l^{\frac{n-2}{2}}}\diff x &=&\wtilde{S}_{HL}\frac{H(a_i,a_i)}{\l^{\frac{n-2}{2}}}\int_{\Real^n}\int_{\Real^n}\frac{U_{(a_i,\l)}^{\soexp{2}}(y) U_{(a_i,\l)}^{\soexp{2}-1}(x)}{\abs{x-y}^\m}\diff x\diff y+O(\frac{1}{\l^n}),\\
&=& \frac{n(n-2)}{\g_0^{\frac{4}{n-2}}}\frac{H(a_i,a_i)}{\l^{\frac{n-2}{2}}}\int_{\Real^n}U_{(a_i,\l)}^{2^\ast-1}\diff x+O(\frac{1}{\l^n}),
\eeqq
since, $-\D U_{(a_i,\l)}=\g_0^{\frac{-4}{n-2}}n(n-2)U_{(a_i,\l)}^{2^\ast-1}$ in $\Real^n$ and therefore
\beq\label{eq2.88}
\int_{\Real^n}\frac{U_{(a_i,\l)}^{\soexp{2}}(y)}{\abs{x-y}^\m}\diff y &=& \frac{n(n-2)}{\wtilde{S}_{HL}\g_0^{\frac{4}{n-2}}}U_{(a_i,\l)}^{2^\ast-\soexp{2}}(x).
\eeq
Using the fact that
\beqq
\int_{\Real^n}U_{(a_i,\l)}^{2^\ast-1}\diff x &=& \frac{\g_0^{2^\ast-1}}{\l^{\frac{n-2}{2}}}\int_{\Real^n}\frac{\diff x}{(1+\abs{x}^2)^{\frac{n+2}{2}}},
\eeqq
we get
\beq\label{eq2.9}
\int_\O-\D U_{(a_i,\l)}\frac{H(a_i,x)}{\l^{\frac{n-2}{2}}}\diff x &=& \g_0 n(n-2)c_1\frac{H(a_i,a_i)}{\l^{n-2}}+O(\frac{1}{\l^n}).
\eeq
From \eqref{eq2.8} and \eqref{eq2.9}, equality \eqref{eq2.7} reduces to
\beq\label{eq2.10}
\abs{PU_{(a_i,\l)}}_{1\O}^2 &=&1- \g_0 n(n-2)c_1\frac{H(a_i,a_i)}{\l^{n-2}}+O(\frac{1}{\l^n}).
\eeq
Now, for $i\neq j$, we have
\beqq
\seq{PU_{(a_i,\l)},PU_{(a_j,\l)}} &=& \int_\O-\D U_{(a_i,\l)}PU_{(a_j,\l)}\diff x\\
&=& \wtilde{S}_{HL}\int_{\O}\int_{\O}\frac{U_{(a_i,\l)}^{\soexp{2}}(y)U_{(a_i,\l)}^{\soexp{2}-1}(x)PU_{(a_j,\l)}(x)}{\abs{x-y}^\m}\diff x\diff y
\eeqq
Using estimates \eqref{eq2.6} and \eqref{eq2.88}, we get
\beqq
\seq{PU_{(a_i,\l)},PU_{(a_j,\l)}} &\hspace{-0.3cm}=& \hspace{-0.2cm}\frac{n(n-2)}{\g_0^{\frac{4}{n-2}}}\left(\int_{\O}U_{(a_i,\l)}^{2^\ast-1}U_{(a_j,\l)}\diff x - \frac{1}{\l^{\frac{n-2}{2}}}\int_{\O}U_{(a_i,\l)}^{2^{\ast}-1}H(a_j,x)\diff x\right)+O(\frac{1}{\l^n}).\\
&\hspace{-2.5cm}=&\hspace{-1cm}\frac{n(n-2)}{\g_0^{\frac{4}{n-2}}}\left(\g_0^{2^\ast}\int_{\O}\d_{(a_i,\l)}^{2^\ast-1}\d_{(a_j,\l)}\diff x - \frac{\g_0^{2^\ast-1}}{\l^{\frac{n-2}{2}}}\int_{\O}\d_{(a_i,\l)}^{2^{\ast}-1}H(a_j,x)\diff x\right) + O(\frac{1}{\l^n}).
\eeqq
From \cite{BL}, we have
\beqq
\int_{\O}\d_{(a_i,\l)}^{2^\ast-1}\d_{(a_j,\l)}\diff x &=& \frac{c_1}{(\l\abs{a_i-a_j})^{n-2}}+o\biggl(\frac{1}{(\l\diff_a)^{n-2}}\biggr),
\eeqq
and by expanding $H(a_j,x)$ around $a_i$, we have
\beqq
\int_{\O}\d_{(a_i,\l)}^{2^\ast-1}H(a_j,x)\diff x &=& c_1\frac{H(a_i,a_j)}{\l^{\frac{n-2}{2}}}+o\biggl(\frac{1}{\l^{\frac{n-2}{2}}}\biggr).
\eeqq
Therefore,
\beq\label{eq2.11}
\seq{PU_{(a_i,\l)},PU_{(a_j,\l)}} &=& n(n-2)\frac{\g_0c_1}{\l^{n-2}}\left( \frac{\g_0}{\abs{a_i-a_j}^{n-2}} - H(a_i,a_j)\right)+o\biggl(\frac{1}{(\l\diff_a)^{n-2}}\biggr),\nonumber\\
&=& n(n-2)\g_0c_1\frac{G(a_i,a_j)}{\l^{n-2}}+o\biggl(\frac{1}{(\l\diff_a)^{n-2}}\biggr).
\eeq
Estimates \eqref{eq2.10} and \eqref{eq2.11} yield
\beqq
N^{\frac{1}{\soexp{2}}} &=& (\sum_{i=1}^{p}\a_i^2)\left(1 + (\sum_{i=1}^{p}\a_i^2)^{-1}n(n-2)\frac{\g_0c_1}{\l^{n-2}}\biggl(\sum_{i\neq j}\a_i\a_jG(a_i,a_j)\right.\\
&-&\left. \sum_{i=1}^{p}\a_i^2H(a_i,a_i)\biggr)\right)+o\biggl(\frac{1}{(\l\diff_a)^{n-2}}\biggr),
\eeqq
and therefore
\beq\label{eq2.12}
N &=& (\sum_{i=1}^{p}\a_i^2)^{\soexp{2}}\left(1 + \soexp{2} (\sum_{i=1}^{p}\a_i^2)^{-1}n(n-2)\frac{\g_0c_1}{\l^{n-2}}\biggl(\sum_{i\neq j}\a_i\a_jG(a_i,a_j)\right.\nonumber\\
&-&\left. \sum_{i=1}^{p}\a_i^2H(a_i,a_i)\biggr)\right)+o\biggl(\frac{1}{(\l\diff_a)^{n-2}}\biggr).
\eeq
We now turn to expand the denominator of $J(\frac{u}{\abs{u}_{1\O}})$. Let $\rh$ be a small positive constant and denote $\ast$ the convolution product. We have
\beq\label{eq2.13}
D &=& \int_\O\int_\O\frac{\bigl(\sum_{i=1}^{p}\a_i P U_{(a_i,\l_i)}(y)\bigr)^{\soexp{2}}\bigl(\sum_{i=1}^{p}\a_i P U_{(a_i,\l_i)}(x)\bigr)^{\soexp{2}}}{\abs{x-y}^\m}\diff x\diff y,\nonumber\\
&=& \int_\O\biggl(\frac{1}{\abs{x}^\m}\ast\bigl(\sum_{i=1}^{p}\a_i P U_{(a_i,\l)}(x)\bigr)^{\soexp{2}}\biggr)\biggl(\sum_{i=1}^{p}\a_i P U_{(a_i,\l)}(x)\biggr)^{\soexp{2}}\diff x,\nonumber\\
&=& \int_{\bigcup_{i=1}^{p}B(a_i,\rh)}\biggl(\frac{1}{\abs{x}^\m}\ast\bigl(\sum_{i=1}^{p}\a_i P U_{(a_i,\l)}(x)\bigr)^{\soexp{2}}\biggr)\biggl(\sum_{i=1}^{p}\a_i P U_{(a_i,\l)}(x)\biggr)^{\soexp{2}}\diff x\nonumber\\
&+& \int_{\O\setminus\bigcup_{i=1}^{p}B(a_i,\rh)}\biggl(\frac{1}{\abs{x}^\m}\ast\bigl(\sum_{i=1}^{p}\a_i P U_{(a_i,\l)}(x)\bigr)^{\soexp{2}}\biggr)\biggl(\sum_{i=1}^{p}\a_i P U_{(a_i,\l)}(x)\biggr)^{\soexp{2}}\diff x
\eeq
On $B(a_i,\rh)$, $i=1,\cdots,p$, we write
\beqq
\sum_{k=1}^{p}\a_iPU_{(a_k,\l)} &=& \a_iU_{(a_i,\l)}+\sum_{k\neq i}\a_kPU_{(a_k,\l)}+\a_i(PU_{(a_i,\l)}-U_{(a_i,\l)}).
\eeqq
Therefore,
\beqq
&& \int_{B(a_i,\rh)}\biggl(\frac{1}{\abs{x}^\m}\ast\bigl(\sum_{k=1}^{p}\a_k P U_{(a_k,\l)}\bigr)^{\soexp{2}}\biggr)\biggl(\sum_{k=1}^{p}\a_k P U_{(a_k,\l)}\biggr)^{\soexp{2}}\diff x \\
&=& \int_{B(a_i,\rh)}\biggl(\frac{1}{\abs{x}^\m}\ast\bigl(\sum_{k=1}^{p}\a_k P U_{(a_i,\l)}\bigr)^{\soexp{2}}\biggr)\bigl(\a_{i}U_{(a_i,\l)}\bigr)^{\soexp{2}}\diff x \\
&+& \soexp{2}\int_{B(a_i,\rh)}\biggl(\frac{1}{\abs{x}^\m}\ast\bigl(\sum_{k=1}^{p}\a_k P U_{(a_k,\l)}\bigr)^{\soexp{2}}\biggr)\bigl(\a_{i}U_{(a_i,\l)}\bigr)^{\soexp{2}-1}\\
&&\biggl(\sum_{k\neq i}\a_kPU_{(a_k,\l)} + \a_i(PU_{(a_i,\l)}-U_{(a_i,\l)})\biggr)\diff x +o\biggl(\frac{1}{(\l\diff_a)^{n-2}}\biggr),
\eeqq
\beqq
&=& \int_{B(a_i,\rh)}\biggl(\frac{1}{\abs{x}^\m}\ast\bigl(\a_i U_{(a_i,\l_i)}\bigr)^{\soexp{2}}\biggr)\bigl(\a_{i}U_{(a_i,\l_i)}\bigr)^{\soexp{2}}\diff x \\
&+& \soexp{2}\int_{B(a_i,\rh)}\biggl(\frac{1}{\abs{x}^\m}\ast\bigl(\a_{i}U_{(a_i,\l)}\bigr)^{\soexp{2}-1}\biggl(\sum_{k\neq i}\a_kPU_{(a_k,\l)}\\
&+& \a_i\bigl(PU_{(a_i,\l)}-U_{(a_i,\l)}\bigr)\biggr)\biggr)\bigl(\a_{i}U_{(a_i,\l)}\bigr)^{\soexp{2}}\diff x \\
&+& \soexp{2}\int_{B(a_i,\rh)}\biggl(\frac{1}{\abs{x}^\m}\ast\bigl(\a_{i}U_{(a_i,\l)}\bigr)^{\soexp{2}}\bigl(\a_{i}U_{(a_i,\l_i)}\bigr)^{\soexp{2} -1}\biggl(\sum_{k\neq i}\a_kPU_{(a_k,\l)} \\
&+& \a_i(PU_{(a_i,\l)}-U_{(a_i,\l)})\biggr)+o\biggl(\frac{1}{(\l\diff_a)^{n-2}}\biggr), \\
&=& \a_i^{2\,\soexp{2}}\int_{\Real^n}\biggl(\frac{1}{\abs{x}^\m}\ast U_{(a_i,\l)}^{\soexp{2}}\biggr)U_{(a_i,\l)}\bigr)^{\soexp{2}}\diff x \\
\eeqq
\beqq
&+& \soexp{2}\int_{\Real^n}\biggl(\frac{1}{\abs{x}^\m}\ast\bigl(\a_{i}U_{(a_i,\l)}\bigr)^{\soexp{2}-1}\biggl(\sum_{k\neq i}\a_kPU_{(a_k,\l)}\biggr)\biggr)\bigl(\a_{i}U_{(a_i,\l)}\bigr)^{\soexp{2}}\diff x\\
&+& \soexp{2}\int_{\Real^n}\biggl(\frac{1}{\abs{x}^\m}\ast\bigl(\a_{i}U_{(a_i,\l)}\bigr)^{\soexp{2}-1}\biggr)\biggl(\a_i\bigl(PU_{(a_i,\l)} - U_{(a_i,\l)}\bigr)\biggr)\bigl(\a_{i}U_{(a_i,\l)}\bigr)^{\soexp{2}}\diff x\bigr)\diff x\\
&+& \soexp{2}\int_{\Real^n}\biggl(\frac{1}{\abs{x}^\m}\ast\bigl(\a_{i}U_{(a_i,\l)}\bigr)^{\soexp{2}}\biggr)\bigl(\a_{i}U_{(a_i,\l)}\bigr)^{\soexp{2} -1}\biggr)\biggl(\sum_{k\neq i}\a_kPU_{(a_k,\l)}\biggr)\diff x\\
&+& \soexp{2}\int_{\Real^n}\biggl(\frac{1}{\abs{x}^\m}\ast\bigl(\a_{i}U_{(a_i,\l)}\bigr)^{\soexp{2}}\biggr)\bigl(\a_{i}U_{(a_i,\l)}\bigr)^{\soexp{2}-1}\biggl(\a_i\bigl(PU_{(a_i,\l)} - U_{(a_i,\l)}\bigr)\biggr)\diff x +o\biggl(\frac{1}{(\l\diff_a)^{n-2}}\biggr).
\eeqq
Using the estimate of \eqref{eq2.6}, we obtain that
\beqq
&& \int_{B(a_i,\rh)}\biggl(\frac{1}{\abs{x}^\m}\ast\bigl(\sum_{k=1}^{p}\a_k P U_{(a_k,\l)}\bigr)^{\soexp{2}}\biggr)\bigl(\sum_{k=1}^{p}\a_k P U_{(a_k,\l)}\bigr)^{\soexp{2}}\diff x \\
&=& \a_i^{2\,\soexp{2}}\wtilde{S}_{HL}^{-1}+\soexp{2}\sum_{k\neq i}\a_i^{2\,\soexp{2}-1}\a_k\int_{\Real^n}\int_{\Real^n}\frac{U_{(a_i,\l)}^{\soexp{2} - 1}(y)PU_{(a_k,\l)}(y)U_{(a_i,\l)}^{\soexp{2}}(x)}{\abs{x-y}^\m}\diff x\diff y\\
&-& \soexp{2}\frac{\a_i^{2\,\soexp{2}}}{\l^{\frac{n-2}{2}}}\int_{\Real^n}\int_{\Real^n}\frac{U_{(a_i,\l)}^{\soexp{2} - 1}(y)H{(a_i,y)}U_{(a_i,\l)}(x)^{\soexp{2}}}{\abs{x-y}^\m}\diff x\diff y\\
&+& \soexp{2}\sum_{k\neq i}\a_i^{2\,\soexp{2}-1}\a_k\int_{\Real^n}\int_{\Real^n}\frac{U_{(a_i,\l)}^{\soexp{2}}(y)U_{(a_i,\l)}^{\soexp{2} - 1}(x)PU_{(a_k,\l)}(x)}{\abs{x-y}^\m}\diff x\diff y\\
&-& \soexp{2}\frac{\a_i^{2\,\soexp{2}}}{\l^{\frac{n-2}{2}}}\int_{\Real^n}\int_{\Real^n}\frac{U_{(a_i,\l)}^{\soexp{2}}(y)U_{(a_i,\l)}^{\soexp{2} -1}(x)H{(a_i,x)}}{\abs{x-y}^\m}\diff x\diff y+ o\biggl(\frac{1}{(\l\diff_a)^{n-2}}\biggr),\\
&=& \a_i^{2\,\soexp{2}}\wtilde{S}_{HL}^{-1} + 2\,\soexp{2}\wtilde{S}_{HL}^{-1}\sum_{k\neq i}\a_i^{2\,\soexp{2} - 1}\a_k\seq{PU_{(a_i,\l)},PU_{(a_k,\l)}}\\
&-& 2\;\soexp{2}\frac{\a_i^{2\,\soexp{2}}}{\l^{\frac{n-2}{2}}}\int_{\Real^n}\int_{\Real^n}\frac{U_{(a_i,\l)}(y)^{\soexp{2}}U_{(a_i,\l)}(x)^{\soexp{2} -1}H{(a_i,x)}}{\abs{x-y}^\m}\diff x\diff y+o\biggl(\frac{1}{(\l\diff_a)^{n-2}}\biggr).
\eeqq
Observe that, by estimate \eqref{eq2.88}, we have
\beqq
\int_{\Real^n}\int_{\Real^n}\frac{U_{(a_i,\l)}^{\soexp{2}}(y)U_{(a_i,\l)}^{\soexp{2} -1}(x)H{(a_i,x)}}{\abs{x-y}^\m}\diff x\diff y
&=& \wtilde{S}_{HL}^{-1}\frac{n(n-2)}{\g_0^{\frac{4}{n-2}}}\int_{\Real^n}U_{(a_i,\l)}^{\soexp{2} -1}(x)H{(a_i,x)}\diff x\\
&=& \wtilde{S}_{HL}^{-1}n(n-2)c_1\g_0\frac{H(a_i,a_i)}{\l^{\frac{n-2}{2}}}+o\biggl(\frac{1}{\l^{\frac{n-2}{2}}}\biggr).
\eeqq
The above estimate with \eqref{eq2.11} yield
\beq\label{eq2.14}
&& \int_{B(a_i,\rh)}\biggl(\frac{1}{\abs{x}^\m}\ast\bigl(\sum_{k=1}^{p}\a_k P U_{(a_k,\l)}\bigr)^{\soexp{2}}\biggr)\bigl(\sum_{k=1}^{p}\a_k P U_{(a_k,\l)}\bigr)^{\soexp{2}}\diff x \\
&=& \wtilde{S}_{HL}^{-1}\a_i^{2\,\soexp{2}-2}\Biggl[\a_i^2 + 2\,\soexp{2}n(n-2)\frac{c_1\g_0}{\l^{n-2}}\biggl(\sum_{k\neq i}\a_i\a_kG(a_i,a_k)-\a_i^2H(a_i,a_i)\biggr)\Biggr] + o\biggl(\frac{1}{(\l\diff_a)^{n-2}}\biggr).\nonumber
\eeq
The remainder integral of \eqref{eq2.13} satisfies
\beq\label{eq2.15}
\hspace{-1cm}\int_{\O\backslash\underset{i=1}{\overset{p}{\bigcup}}B(a_i,\rh)}\bigl(\frac{1}{\abs{x}^\m}\ast\bigl(\sum_{k=1}^{p}\a_i P U_{(a_i,\l)}(x)\bigr)^{\soexp{2}}\bigr)\bigl(\sum_{k=1}^{p}\a_i P U_{(a_i,\l)}(x)\bigr)^{\soexp{2}}\diff x&=& O\bigl(\frac{1}{(\l\diff_a)^{n}}\bigr).
\eeq
Using \eqref{eq2.14} and \eqref{eq2.15}, the expansion of \eqref{eq2.13} reduces to
\beq\label{eq2.16}
D &=& \wtilde{S}_{HL}^{-1}\bigl(\sum_{i=1}^{p}\a_i^{2\,\soexp{2}}\bigr)\Biggl[1 + 2\,\soexp{2}n(n-2)\frac{c_1\g_0}{\l^{n-2}}\bigl(\sum_{k\neq i}\frac{\a_i^{2\,\soexp{2}-1}\a_k}{\sum_{j=1}^{p}\a_j^{2\,\soexp{2}}}G(a_i,a_k) \nonumber\\
&-& \sum_{i=1}^{p}\frac{\a_i^{2\,\soexp{2}}}{\sum_{j=1}^{p}\a_j^{2\,\soexp{2}}}H(a_i,a_i)\bigr)\Biggr] + o\bigl(\frac{1}{(\l\diff_a)^{n-2}}\bigr).
\eeq
The expansion of Proposition \ref{prop2.5} follows from \eqref{eq2.12} and \eqref{eq2.16}.
\end{pf}

\n Of course the above expansion is useful when $\diff_a=\underset{i\neq j}{\min}\abs{a_i-a_j}$ is lower bounded by a positive constant independent of $a$. In this case, for any $\e >0$, there exists $\overline{\l}=\overline{\l}(\e)>0$ (independent of $a$, as long as $a$ lies in a compact set of $\O$) such that, for any $\l\geq\overline{\l}$,
\beq\label{eq2.17}
J\biggl(\frac{\sum_{i=1}^{p}\a_iPU_{(a_i,\l)}}{\abs{\sum_{i=1}^{p}\a_iPU_{(a_i,\l)}}_{1\O}}\biggr) &\leq& (p+\e)^{\soexp{2}-1}\wtilde{S}_{HL},
\eeq
since $\frac{(\underset{i=1}{\overset{p}{\sum}}\a_i^2)^{\soexp{2}}}{\underset{i=1}{\overset{p}{\sum}}\a_i^{2\,\soexp{2}}}\leq p^{\soexp{2} - 1}$.

\n If $\diff_a$ is close to zero, or if a parameter $\a_i$ is close to zero for some index $i,$ $i=1,\cdots,p$, we have
\beq\label{eq2.18}
J\biggl(\frac{\sum_{i=1}^{p}\a_iPU_{(a_i,\l)}}{\abs{\sum_{i=1}^{p}\a_iPU_{(a_i,\l)}}_{1\O}}\biggr) &\leq& p^{\soexp{2}-1}\wtilde{S}_{HL}.
\eeq
This can be proved as the analogous assertion in \cite{BC 1989}.

\section{Proof of the existence Theorem}\label{sec3}
\n To prove Theorem \ref{th1.1} we argue by contradiction. Therefore we suppose in what follows that problem \eqref{eq1.3} has no positive solution.

\n Let $-\pa J\colon\Sig\to T\Sig$ be the gradient of $J$. It is a continuous vector field, locally Lipschitz and bounded. For $u\in\Sig$, we denote by $s\longmapsto \eta(s,u)$, $s\geq 0$, the unique solution of the differential equation
\beqq
\left\{
  \begin{array}{ll}
    \dot{\eta} &=-\pa J(\eta) \hbox{} \\
    \eta(0) &=u \hbox{.}
  \end{array}
\right.
\eeqq
It is straightforward to see that $\eta(s,u)$ is defined for any $s\in [0,\infty)$, $J(\eta(s,u))$ is a decreasing function and $\pa J(\eta(s,u))\to 0$, as $s\to\infty$. Moreover, if we denote by $\Sig^+=\set{u\in\Sig,\;u\geq 0}$, then $\Sig^+$ is invariant under the action of $(-\pa J)$, in the sense that, if $u\in\Sig^+$ then $\eta(s,u)\in\Sig^+$ for any $s\geq 0$. For any $p\in\Natl\backslash\{0\}$ and $\e>0$, we define
\beqq
V(p,\e) &=& \Bigg\{u\in\Sig^+,\;\exists a_1,\cdots,a_p\in\O,\;\exists\l_1,\cdots,\l_p>0\text{  and  }\a_1,\cdots,\a_p>0,\;s.t.\,, \\
&&|u - \sum_{i=1}^{p}\a_iPU_{(a_i,\l_i)}|_{1\O}<\e\text{  with  }\l_i>\e^{-1},\;\l_i\diff(a_i,\pa\O)>\e^{-1},\;\forall i=1,\cdots,p\\
&&\text{and  }\e_{ij}=\bigl(\frac{\l_i}{\l_j} + \frac{\l_j}{\l_i} + \l_i\abs{a_i - a_j}^2\bigr)^{\frac{2-n}{2}}<\e,\;\forall i\neq j\Bigg\}.
\eeqq
We have the following result
\begin{proposition}\label{prop3.1}
Assume that $J$ has no positive critical point. Let $u\in\Sig^+$. There exists a unique positive integer $p$ which depends on $u$ only, such that for any $\e>0$ there exists $s_\e>0$ such that for any $s>s_\e$, $\eta(s,u)\in V(p,\e)$. Moreover $J(\eta(s,u))\to\ell$, as $s\to\infty$, with $\ell\leq p^{\soexp{2}-1}\wtilde{S}_{HL}$.
\end{proposition}
\begin{pf}
Let $(s_k)_k$ be a positive sequence tending to $\infty$. Denote $u_k=\eta(s_k,u)$, $u\in\Sig^+$. Therefore, $(u_k)_k$ is a sequence of $\Sig^+$ which satisfies
\beqq
J(u_k)\to c &\text{ and }& \pa J(u_k)\to 0.
\eeqq
Since we have supposed that $J$ has no critical point in $\Sig^+$, by (\cite{BC 1991}, Proposition 1) and (\cite{31 Aut}, Lemma 4.5), Here exists $p\in\Natl\backslash\{0\}$ such that $u_k\in V(p,\e_k)$, $\forall k$, where $\e_k>0$ and $\e_k\to 0$ as $k\to\infty$. Fix $\e>0$. If the flow line $\eta(s,u)$ enter $V(p,\frac{\e}{2})$ at a time $s_k$ and leave $V(p,\e)$ at a later time $t_k$, for $t\in[s_k,t_k]$, $\eta(s,u)\in V(p,\e)\backslash V(p,\frac{\e}{2})$. Therefore,
\be\label{eq3.1}
J(\eta(t_k,u))-J(\eta(s_k,u))=\int_{s_k}^{t_k}\dot{\widehat{J(\eta(s,u))}}\diff s\leq -\int_{s_k}^{t_k}\abs{\pa J(\eta(s,u))}_{1\O}^2\diff s\leq -\g(t_k-s_k),
\ee
where $\g=\g(\e)$ is the lower bound of $\abs{\pa J(u)}_{1\O}$ on $V(p,\e)\backslash V(p,\frac{\e}{2})$. By (\cite{BC 1991}, Proposition 1) and (\cite{31 Aut}, Lemma 4.5), $\g$ is positive. Now, let
\beqq
0<\b &=& \diff(V(p,\e),V(p,\e)^c).
\eeqq
We have
\beqq
\b\leq\diff(\eta(s_k,u),\eta(t_k,u)) &\leq& \int_{s_k}^{t_k}\abs{\dot{\eta}(s,u)}_{1\O}\diff s\leq M_0(t_k-s_k),
\eeqq
where $M_0=\underset{u\in V(p,\e)}{\sup}\abs{\pa J(u)}$ and hence $t_k-s_k\geq\frac{\b}{M_0}$. This with \eqref{eq3.1} yield
\be\label{eq3.2}
J(\eta(t_k,u))-J(\eta(s_k,u)) \leq -\frac{\g\b}{M_0}
\ee
It follows from \eqref{eq3.2} that there is a finite number of such intervals $[s_k,t_k]$, since $J$ is lower bounded. As a consequence, there exists $s_\e>0$ such that $\eta(s,u)\in V(p,\e)$, $\forall s\geq s_\e$. The limit of $J(\eta(s,u))$ follows from an expansion like the one of Proposition \ref{prop2.5}.
\end{pf}

\n In the following we denote by $J_c$, $c\in\Real$, the set
\[J_c = \set{u\in\Sig^+,\;J(u)\leq c}.\]
Let $p\in\Natl\backslash\{0\}$ and let $\e_p$ be a constant such that
\[0<\e_p<\wtilde{S}_{HL}((p+1)^{\soexp{2}-1}-p^{\soexp{2}-1}).\]
Denote,
\[W_u(V(p,\e_p)) = \set{\eta(s,u),\;s\geq 0,\;u\in\bigcup_{1\leq k\leq p} V(k,\e_k)}.\]
We define,
\beqq
E_p=J_{p^{\soexp{2}-1}\wtilde{S}_{HL}+\e_p}\cap W_u(V(p,\e_p)) &\text{ and }& F_p=J_{p^{\soexp{2}-1}\wtilde{S}_{HL}}\cap W_u(V(p,\e_p)).
\eeqq
In this way,
\[F_{p-1}\subset E_{p-1}\subset F_p\subset E_p,\;\;\forall p\geq 1,\]
with $F_0=E_0=\emptyset$.

\n For a pair of topological $(X,Y)$, $Y\subset X$, we denote $H_\ell(X,Y)$, $\ell\geq 0$, the relative homology of $(X,Y)$. We then have
\begin{proposition}\label{prop3.2}
Let $p\in\Natl\backslash\{0\}$. For any $\e\in(0,\e_p)$, there exists a continuous mapping
\beqq
R\colon\;(E_p,F_p) &\to& \biggl(J_{p^{\soexp{2}-1}\wtilde{S}_{HL}+\e}\cap V(p,\e_p),J_{p^{\soexp{2}-1}\wtilde{S}_{HL}}\cap V(p,\e_p)\biggr),
\eeqq
which induces a sequence of isomorphisms
\beqq
R_\ell\colon\;H_\ell(E_p,F_p) &\to& H_\ell\biggl(J_{p^{\soexp{2}-1}\wtilde{S}_{HL}+\e}\cap V(p,\e_p),J_{p^{\soexp{2}-1}\wtilde{S}_{HL}}\cap V(p,\e_p)\biggr),\;\;\ell\geq 0.
\eeqq
\end{proposition}
\begin{pf}
Let $0<\e<\e_p$. Using the gradient flow $\eta(.,.)$ of the field $(-\pa J)$, the following retract by deformation holds,
\beqq
J_{p^{\soexp{2}-1}\wtilde{S}_{HL}+\e_p} &\simeq& J_{p^{\soexp{2}-1}\wtilde{S}_{HL}+\e},
\eeqq
since by Proposition 4.1, $J$ has no critical value at infinity between the levels $p^{\soexp{2}-1}\wtilde{S}_{HL}+\e_p$ and $p^{\soexp{2}-1}\wtilde{S}_{HL}+\e$. Here $\simeq$ denotes retract by deformation. It follows that
\beqq
(E_p,F_p) &\simeq& \biggl(J_{p^{\soexp{2}-1}\wtilde{S}_{HL}+\e}\cap W_u(V(p,\e_p)),F_p\biggr).
\eeqq
Let us denote by
\beqq
r\colon\;(E_p,F_p) &\to& \biggl(J_{p^{\soexp{2}-1}\wtilde{S}_{HL}+\e}\cap W_u(V(p,\e_p)),F_p\biggr)
\eeqq
the corresponding retract by deformation. Using the fact that
\beqq
J_{p^{\soexp{2}-1}\wtilde{S}_{HL}+\e}\cap W_u(V(p,\e_p))\backslash F_p &\subset& V(p,\e_p),
\eeqq
we obtain that
\beqq
J_{p^{\soexp{2}-1}\wtilde{S}_{HL}+\e}\cap W_u(V(p,\e_p))\backslash F_p &=& J_{p^{\soexp{2}-1}\wtilde{S}_{HL}+\e}\cap V(p,\e_p)\backslash J_{p^{\soexp{2}-1}\wtilde{S}_{HL}}\cap V(p,\e_p),
\eeqq
and hence the following two pairs
\beqq
\biggl(J_{p^{\soexp{2}-1}\wtilde{S}_{HL}+\e}\cap W_u(V(p,\e_p)),F_p\biggr) &\text{ and }& \biggl(J_{p^{\soexp{2}-1}\wtilde{S}_{HL}+\e}\cap V(p,\e_p),J_{p^{\soexp{2}-1}\wtilde{S}_{HL}}\cap V(p,\e_p)\biggr),
\eeqq
are homotopically equivalent. Denote by
\beqq
\widehat{r}\colon\;\biggl(J_{p^{\soexp{2}-1}\wtilde{S}_{HL}+\e}\cap W_u(V(p,\e_p)),F_p\biggr) &\to& \biggl(J_{p^{\soexp{2}-1}\wtilde{S}_{HL}+\e}\cap V(p,\e_p),J_{p^{\soexp{2}-1}\wtilde{S}_{HL}}\cap V(p,\e_p)\biggr),
\eeqq
the corresponding homotopy equivalence. Setting
\beqq
R &=& \widehat{r}o r.
\eeqq
$R$ is an homotopy equivalence, since $r$ and $\widehat{r}$ are. Hence the result follows.
\end{pf}

\n We now introduce the following notations. For $p\in\Natl\backslash\{0\}$, we denote by
\beqq
\D_{p-1} &=& \Bigl\{(\a_1,\cdots,\a_p)\in[0,1]^p,\;\;\sum_{i=1}^{p}\a_i=1\Bigr\}
\eeqq
the standard $p$-simplex in $\Real^p$. Let $M$ be a fixed compact $k$-dimensional manifold in $\O$, $1\leq k\leq n-1$. We denote
\beqq
B_p(M) &=& \set{\sum_{i=1}^{p}\a_i\d_{a_i},\;\a_i\in[0,1],\;a_i\in M,\;\forall i=1,\cdots,p,\;\text{ and  }\sum_{i=1}^{p}\a_i=1},
\eeqq
where $\d_{a_i}$ denotes the Dirac mass at $a_i$. We agree that $B_0(M)=\emptyset$. Lastly, denote $\varrho_p$ the symmetric group of order $p$. We then have.
\begin{proposition}\label{prop3.3}
For any $P\in\Natl\backslash\{0\}$, the homology $H_\ell(E_p,F_p)$, $\ell\geq 0$, of the pair $(E_p,F_p)$, is a module on the cohomology $H^\ell(\O^p/\varrho_p)$, $\ell\geq 0$ of $\O^p/\varrho_p$. Moreover there exits a continuous mapping of topological pairs
\beqq
\phi_p\colon (B_p(M),B_{p-1}(M)) &\to& (E_p,F_p),
\eeqq
such that the homology homorphism
\beqq
(\phi_p)_\ell\colon H_\ell(B_p(M),B_{p-1}(M)) &\to& (E_p,F_p),\;\;\ell\geq 0,
\eeqq
is $H^\ell(\O^p/\varrho_p)$-linear.
\end{proposition}
\begin{pf}
Let $p\in\Natl\backslash\{0\}$. For $\e_p$ small enough and for $u\in V(p,\e_p)$, we optimize the approximation of $u$ with respect to $\underset{i=1}{\overset{p}{\sum}}\a_iPU_{(a_i,\l_i)}$. Namely, we set the following minimization problem
\beqq
\min_{\a_i,a_i,\l_i}\abs{u-\underset{i=1}{\overset{p}{\sum}}\a_iPU_{(a_i,\l_i)}}_{1\O}.
\eeqq
It is proved in \cite{BL} that this problem has a unique solution $\a_i,a_i,\l_i$, $i=1,\cdots,p$ (modulo a permutation). It follows that there exists a continuous mapping
\beq\label{eq3.3}
\chi\colon\;V(p,\e_p) &\to& \O^p/\varrho_p,
\eeq
such that to any $u\in V(p,\e_p)$, we associate $(a_1,\cdots,a_p)$; the unique solution of the above problem of minimization.

\n Let $0<\e<\e_p$ and $\chi_1$ be the restriction of $\chi$ on $J_{p^{\soexp{2}-1}\wtilde{S}_{HL}+\e}\cap V(p,\e_p)$. $\chi_1$ is continuous and induces a cohomology homomorphism
\beqq
(\chi_1)^\ell\colon\;H^\ell(\O^p/\varrho_p) &\to& H^\ell\bigl(J_{p^{\soexp{2}-1}\wtilde{S}_{HL}+\e}\cap V(p,\e_p)\bigr),\;\;\ell\geq 0.
\eeqq
From \cite{2 BB}, we know that for any pair of topological spaces $(X,Y)$, the group of cohomology $H^\ell(X)$, $\ell\geq 0$, of $X$ acts through the cap-product on the relative homology group $H_\ell(X,Y)$ of the pair $(X,Y)$ and the action is linear. We apply this in our statement. We derive that
$$H_\ell\biggl(J_{p^{\soexp{2}-1}\wtilde{S}_{HL}+\e}\cap V(p,\e_p),J_{p^{\soexp{2}-1}\wtilde{S}_{HL}}\cap V(p,\e_p)\biggr),\;\;\ell\geq 0,$$ is a module over
\n $H_\ell\biggl(J_{p^{\soexp{2}-1}\wtilde{S}_{HL}+\e}\cap V(p,\e_p)\biggr)$ and hence over $H^\ell(\O^p/\varrho_p)$ by $(\chi_1)^\ell$. Using the isomorphism $R_\ell$ of Proposition \ref{prop3.2}, we deduce that $H_\ell(E_p,F_p)$ is a module over $H^\ell(\O^p/\varrho_p)$.

\n To construct the homology homomorphism $(\phi_p)_\ell$, $\ell\geq 0$, we first define on $M^p\times\D_{p-1}$ an equivalence relation such that the class of an element $(a_1,\cdots,a_p,\a_i,\cdots,\a_p)\in M^p\times\D_{p-1}$ equals to
\[\set{(a_{\s(1)},\cdots,a_{\s(p)},\a_{\s(1)},\cdots,\a_{\s(p)}),\;\s\in\varrho_p}.\]
Denote $M^p\underset{\varrho_p}{\times}\D_{p-1}$, the related quotient space. We define
\beqq
\left.
  \begin{array}{ll}
    \pi_p\colon\;B_p(M) &\to M^p\underset{\varrho_p}{\times}\D_{p-1} \hbox{} \\
     & \hbox{} \\
    \sum_{i=1}^{p}\a_i\d_i &\longmapsto (a_1,\cdots,a_p,\a_i,\cdots,\a_p) \hbox{.}
  \end{array}
\right.
\eeqq
It is easy to verify that $\pi_p$ is an homeomorphism. Let
\beqq
S_p(M) &=& \set{(a_1,\cdots,a_p)\in M^p,\;\;\exists i\neq j,\;,s.t.,\;a_i=a_j}.
\eeqq
$\pi_p$ induces an homeomorphism of pairs that we denote again $\pi_p$,
\beqq
\pi_p\colon\;(B_p(M),B_{p-1}(M)) &\to& (M^p\underset{\varrho_p}{\times}\D_{p-1},S_p(M)\times\D_{p-1}\underset{\varrho_p}{\cup}M^p\times\pa\D_{p-1}),
\eeqq
where $\pa\D_{p-1}$ denote the boundary of $\D_{p-1}$. Let
\beqq
(\pi_p)_\ell\colon\;H_\ell(B_p(M),B_{p-1}(M)) &\to& H_\ell (M^p\underset{\varrho_p}{\times}\D_{p-1},S_p(M)\times\D_{p-1}\underset{\varrho_p}{\cup}M^p\times\pa\D_{p-1}),\;\ell\geq 0
\eeqq
be the homology isomorphism induced by $\pi_p$. Denote $T_p$ a small $\varrho_p$-equivariant tubular open neighborhood of $S_p(M)$ in $M^p$ such that $S_p(M)$ is a retract by deformation of $T_p$ (see\cite{2 BB}) and denote
\beqq
i_p\colon\;(M^p\underset{\varrho_p}{\times}\D_{p-1},S_p(M)\times\D_{p-1}\underset{\varrho_p}{\cup}M^p\times\pa\D_{p-1}) &\hspace{-0.3cm}\to& \hspace{-0.2cm}(M^p\underset{\varrho_p}{\times}\D_{p-1},\overline{T_p}\times\D_{p-1}\underset{\varrho_p}{\cup}M^p\times\pa\D_{p-1}),
\eeqq
the natural injection. It is an homotopy equivalence, since $\overline{T_p}$ retracts by deformation on $S_p(M)$. Therefore induces an homology isomorphism
\beqq
(i_p)_\ell\colon\;H_\ell(M^p\underset{\varrho_p}{\times}\D_{p-1},S_p(M)\times\D_{p-1}\underset{\varrho_p}{\cup}M^p\times\pa\D_{p-1}) &\hspace{-0.3cm}\to& \hspace{-0.2cm} H_\ell(M^p\underset{\varrho_p}{\times}\D_{p-1},\overline{T_p}\times\D_{p-1}\underset{\varrho_p}{\cup}M^p\times\pa\D_{p-1}).
\eeqq
Let $\th$ be a given small positive constant, We denote
\beqq
\D_{p-1}^\th &=& \set{(\a_i,\cdots,\a_p)\in\D_{p-1},\frac{\a_i}{\a_j}\in[1-\th,1+\th],\;\forall i\neq j}.
\eeqq
Using the fact that $(\D_{p-1}^\th)^c$ retracts on $\pa\D_{p-1}$, there exists an homology isomorphism
\beqq
(j_p)_\ell\colon\;H_\ell(M^p\underset{\varrho_p}{\times}\D_{p-1},\ov{T_p}\times\D_{p-1}\underset{\varrho_p}{\cup}M^p\times\pa\D_{p-1}) &\hspace{-0.3cm}\to& \hspace{-0.2cm} H_\ell(M^p\underset{\varrho_p}{\times}\D_{p-1},\overline{T_p}\times\D_{p-1}\underset{\varrho_p}{\cup}M^p\times(\D_{p-1}^\th)^c),
\eeqq
$\ell\geq 0$. Let us denote
\beq\label{eq3.33}
M_0^p &=& M^p\backslash T_p
\eeq
By excision \cite{4 BB} of $T_p\underset{\varrho_p}{\times}\D_{p-1}\underset{\varrho_p}{\cup}T_p\times(\D_{p-1}^\th)^c$, the two pairs

\n$(M^p\underset{\varrho_p}{\times}\D_{p-1},\overline{T_p}\times\D_{p-1}\underset{\varrho_p}{\cup}M^p\times(\pa\D_{p-1}^\th)^c)$ and $(M_0^p\underset{\varrho_p}{\times}\D_{p-1},\pa M_0^p\times\D_{p-1}\underset{\varrho_p}{\cup}M_0^p\times(\D_{p-1}^\th)^c)$

\n have the same type of homotopy. Therefore

\n there exists a mapping
\[h_p\colon\,(M^p\underset{\varrho_p}{\times}\D_{p-1},\overline{T_p}\times\D_{p-1}\underset{\varrho_p}{\cup}M^p\times(\D_{p-1}^\th)^c) \to (M_0^p\underset{\varrho_p}{\times}\D_{p-1},\pa M_0^p\times\D_{p-1}\underset{\varrho_p}{\cup}M_0^p\times(\D_{p-1}^\th)^c)\]
which defines an homotopy equivalence. Therefore, the corresponding homology homomorphism $(h_p)_\ell$, $\ell\geq 0$, is indeed an isomorphism.

\n Let for $\l>0$, $f_{p,\l}$ denotes
\beqq
&&\left.
  \begin{array}{ll}
    \hspace{1.cm}f_{p,\l}\colon M^p\underset{\varrho_p}{\times}\D_{p-1} &\longrightarrow \Sig^+ \hbox{} \\
     & \hbox{} \\
    (a_1,\cdots,a_p,\a_i,\cdots,\a_p) &\longmapsto \frac{\sum_{i=1}^{p}\a_iPU_{(a_i,\l)}}{\abs{\sum_{i=1}^{p}\a_iPU_{(a_i,\l)}}_{1\O}} \hbox{.}
  \end{array}
\right.
 \\
\eeqq
Using the fact $\diff(a_i,a_j)$, $1\leq i\neq j\leq p$ lower bounded by $\diff_1>0$, for every $(a_1,\cdots,a_p)\in M_0^p$, where $\diff_1$ is the diameter of $T_p$, we can select $\ov{\l}=\ov{\l}(p)$ large enough, so that for any $\l\geq\ov{\l}$,
\beq\label{eq3.4}
f_{p,\l}( M_0^p\underset{\varrho_p}{\times}\D_{p-1})\subset &&\hspace{-0.7cm}\underset{1\leq k\leq p}{\cup}V(k,\e_k)\subset W_u(V(p,\e_p)),
\eeq
and
\beq\label{eq3.5}
f_{\l}(\pa M_0^p\underset{\varrho_p}{\times}\D_{p-1}\underset{\varrho_p}{\cup}M_0^p\times(\D_{p-1}^\th)^c)\subset &&\hspace{-0.7cm}\underset{1\leq k\leq p}{\cup}V(k,\e_k)\subset W_u(V(p,\e_p)),
\eeq
Moreover, by estimate \eqref{eq2.17}, we have
\beq\label{eq3.6}
f_{p,\l}( M_0^p\times\D_{p-1})\subset &&\hspace{-0.7cm}J_{p^{\soexp{2}-1}\wtilde{S}_{HL}+\e_p},
\eeq
and by estimate \eqref{eq2.18}, we have
\beq\label{eq3.7}
f_{p,\l}(\pa M_0^p\times\D_{p-1}\underset{\varrho_p}{\cup}M_0^p\times(\pa\D_{p-1}^\th)^c)\subset &&\hspace{-0.7cm}J_{p^{\soexp{2}-1}\wtilde{S}_{HL}},
\eeq
provided $\th$ and $\diff_1$ are small. It follows from \eqref{eq3.4}--\eqref{eq3.7} that $f_{p,\l}$ maps the pair $$(M_0^p\underset{\varrho_p}{\times}\D_{p-1},\pa M_0^p\underset{\varrho_p}{\times}\D_{p-1}\underset{\varrho_p}{\cup}M_0^p\times(\D_{p-1}^\th)^c)$$
into $(E_p,F_p)$, we denote again
\beqq
f_{p,\l}\colon\,(M_0^p\underset{\varrho_p}{\times}\D_{p-1},\pa M_0^p\times\D_{p-1}\underset{\varrho_p}{\cup}M_0^p\times(\D_{p-1}^\th)^c) &\to& (E_p,F_p)
\eeqq
the corresponding mapping. Passing to homology, we obtain the following homomorphism
\beqq
(f_{p,\l})_\ell\colon\,H_\ell(M_0^p\underset{\varrho_p}{\times}\D_{p-1},\pa M_0^p\times\D_{p-1}\underset{\varrho_p}{\cup}M_0^p\times(\pa\D_{p-1}^\th)^c) &\to& (E_p,F_p),\;\ell\geq 0.
\eeqq
We now define the required mapping $\phi_p$ by:
\[\phi_p=f_{p,\l}oh_poj_poi_po\pi_p.\]
By construction $\phi_p$ is continuous and hence induces an homology homomorphism $(\phi_p)_\ell$, $\ell\geq 0$, from $H_\ell(B_p(M),B_{p-1}(M))$ into $H_\ell(E_p,F_p)$.

\n Let us observe that $f_{p,\l}$ maps for $\th>0$ small enough and $\l>\ov{\l}$, the sets
\beq\label{eq3.8}
M_0^p\underset{\varrho_p}{\times}\D_{p-1} &\text{ into }& E_p\cap V(p,\e_p),
\eeq
and
\beq\label{eq3.9}
\pa(M_0^p\underset{\varrho_p}{\times}\D_{p-1})= \pa M_0^p\times\D_{p-1}^\th\underset{\varrho_p}{\cup}M_0^p\times\pa\D_{p-1}^\th &\text{ into }& F_p\cap V(p,\e_p).
\eeq
Moreover, if we consider $M^p\underset{\varrho_p}{\times}\D_{p-1}$ as $B_p(M)$, through the homeomorphism $\pi_p$, we then have
\[(B_p(M),(M_0^p\underset{\varrho_p}{\times}\D_{p-1}^\th)^c)\simeq(B_p(M),B_{p-1}(M)),\]
since $(M_0^p)^c=\ov{T_p}$ retracts by deformation on $S_p(M)$ and $(\D_{p-1}^\th)^c$ retracts on $\pa\D_{p-1}$. By excision

\vspace{0.5cm}\n of $\tikzarc{(M_0^p\underset{\varrho}{\times}\D_{p-1}^\th)^c}$, the two pairs $(M_0^p\underset{\varrho_p}{\times}\D_{p-1}^\th,\pa(M_0^p\underset{\varrho_p}{\times}\D_{p-1}^\th))$ and $(B_p(M),B_{p-1}(M))$ have the same type of homology. Let
\beqq
q_p\colon\,(M_0^p\underset{\varrho_p}{\times}\D_{p-1}^\th,\pa(M_0^p\underset{\varrho_p}{\times}\D_{p-1}^\th)) &\to& (B_p(M),B_{p-1}(M))
\eeqq
be the associated homotopy equivalence.

\n Consider now the diagram of continuous maps,
\[\xymatrix@C=2.5cm@R=1.cm{
(B_p(M),B_{p-1}(M)) \ar[r]^{\phi_p}  & (E_p,F_p)  \\
(M_0^p\underset{\varrho_p}{\times}\D_{p-1}^\th,\pa(M_0^p\underset{\varrho_p}{\times}\D_{p-1}^\th)) \ar[r] ^{f_{p,\l}} \ar[u]_{q_p} \ar[d]^{\widehat{i}_p} & (E_p\cap V(p,\e_p),F_p\cap V(p,\e_p)) \ar[u]_{\wtilde{i}_p} \ar[d]^{\chi}\\
M^p\underset{\varrho_p}{\times}\D_{p-1} \ar[r]^{V_p} & \O^p/\varrho_p,
}\]
where $V_p$ denotes the first projection, $\wtilde{i}_p$ and $\widehat{i}_p$ are the natural injection, $\chi$ is defined in \eqref{eq3.3}, $f_{p,\l}$ is defined in \eqref{eq3.8} and \eqref{eq3.9} and $q_p$ and $\phi_p$ are defined in the preceding steps.

\n Passing the above diagram to homology. Using the fact that $(\wtilde{i}_p)_\ell$, $\ell\geq 0$, is an isomorphism and the homology isomorphism $(q_p)_\ell$ can be expressed by
\beq\label{eq3.10}
(q_p)_\ell &=& (i)_\ell o(h_p)_\ell^{-1}o(j_p)_\ell^{-1}o(i_p)_\ell^{-1}o(\pi_p)_\ell^{-1},\;\ell\geq 0,
\eeq
where
\beqq
i\colon\,(M_0^p\underset{\varrho_p}{\times}\D_{p-1}^\th,\pa(M_0^p\underset{\varrho_p}{\times}\D_{p-1}^\th)) &\to& (M_0^p\underset{\varrho_p}{\times}\D_{p-1},\pa M_0^p\times\D_{p-1}\underset{\varrho_p}{\cup}M_0^p\times(\D_{p-1}^\th)^c),
\eeqq
denotes the natural injection, the obtained homological diagram is commutative and hence the homology homomorphism $(\phi_p)_\ell$ , $\ell\geq 0$, is $H^\ell(\O/_{\varrho_p})$-linear, in the sense that for any $\varphi^\ell\in H^\ell(\O^p/_{\varrho_p})$, and $[N_p,N_{p-1}]_\ell\in H_\ell(B_p(M),B_{p-1}(M))$, it holds
\beq\label{eq3.11}
(\phi_p)_\ell(\varphi^\ell.[N_p,N_{p-1}]_\ell) &=& \varphi^\ell.(\phi_p)_\ell([N_p,N_{p-1}]_\ell).
\eeq
\n Note that, it is proved that $H_\ell(E_p,F_p)$ is a module over $H^\ell(\O^p/_{\varrho_p})$ and by the cap-product, $H_\ell(B_p(M),B_{p-1}(M))$ is a module over $H^\ell(B_p(M))$ and thus over $H^\ell(\O^p/_{\varrho_p})$ via the first projection. the proof of the Proposition is thereby completed.
\end{pf}

\n Let us observe that from \eqref{eq3.10}, $(q_p)_\ell$ defines an isomorphism between
\beqq
H_\ell(B_p(M),B_{p-1}(M)) &\text{ and }& H_\ell(M_0^p\underset{\varrho_p}{\times}\D_{p-1}^\th,\pa(M_0^p\underset{\varrho_p}{\times}\D_{p-1}^\th)),\;\ell\geq 0.
\eeqq
Here $M_0^p$ is defined in \eqref{eq3.33}. We also note that $M_0^p\underset{\varrho_p}{\times}\D_{p-1}$ is a manifold of dimension $kp + p-1$, $(k=dim\,M)$. It can be viewed as a singular $(k_p + p-1)$-chain with $\Z_2$ coefficients. Therefore the pair $(M_0^p\underset{\varrho_p}{\times}\D_{p-1},\pa(M_0^p\underset{\varrho_p}{\times}\D_{p-1}))$ defines a $(k_p + p-1)$-cycle in its self and thus gives arise to a non zero class denoted $[M_0^p\underset{\varrho_p}{\times}\D_{p-1},\pa(M_0^p\underset{\varrho_p}{\times}\D_{p-1})]$ in the relative homology of the pair $(M_0^p\underset{\varrho_p}{\times}\D_{p-1},\pa(M_0^p\underset{\varrho_p}{\times}\D_{p-1}))$. In the next, we denote
\beq\label{eq3.12}
[B_p(M),B_{p-1}(M)] &=& (q_p)_\ell([M_0^p\underset{\varrho_p}{\times}\D_{p-1},\pa(M_0^p\underset{\varrho_p}{\times}\D_{p-1})]).
\eeq
Since $(q_p)_\ell$ is an isomorphism, then $[B_p(M),B_{p-1}(M)]$ defines a non zero class in the relative homology group $H_\ell(B_p(M),B_{p-1}(M))$ of the pair $(B_p(M),B_{p-1}(M))$, $\forall p\in\Natl\backslash\{0\}$. We have the following result.
\begin{proposition}
    Let $\phi_p$, $p\in\Natl\backslash\{0\}$, be the continuous mapping defined in Proposition \ref{prop3.3} and let $(\phi_p)_\ell$, $\ell\geq 0$, be the associated homology homomorphism. We then have
\beqq
(\phi_1)_\ell([B_1(M),B_{0}(M)])\neq 0 &\Rightarrow& (\phi_p)_\ell([B_p(M),B_{p-1}(M)])\neq 0\;\forall p\in\Natl\backslash\{0\},
\eeqq
where $[B_p(M),B_{p-1}(M)]$ is defined in \eqref{eq3.12}.
\end{proposition}
\begin{pf}
We argue by induction. For $p=1$, we have
\beqq
(\phi_1)_\ell([B_1(M),B_{0}(M)])&\neq& 0.
\eeqq
Let $p\in\Natl\backslash\{0\}$. Assume that
\beq\label{eq3.13}
(\phi_p)_\ell([B_p(M),B_{p-1}(M)])&\neq& 0.
\eeq
We claim that
\beq\label{eq3.14}
(\phi_{p+1})_\ell([B_{p+1}(M),B_{p}(M)])&\neq& 0.
\eeq
For this, let
\beqq
\pa_{p+1}\colon\,H_{\ell+1}(B_{p+1}(M),B_{p}(M)) &\to& H_\ell(B_p(M),B_{p-1}(M)),\;\ell\geq 0,
\eeqq
be the connecting homomorphism, see \cite{4 BB}, and let $\varphi_{p+1}\in H^\ell(\O^{p+1}/_{\varrho_{p+1}})$ such that
\beqq
\pa_{p+1}(\varphi_{p+1}.[B_{p+1}(M),B_{p}(M)]) &=& [B_{p}(M),B_{p-1}(M)].
\eeqq
See (\cite{BC 1989},Appendix C) for the existence of $\varphi_{p+1}$. By \eqref{eq3.13}, we have
\beq\label{eq3.15}
(\phi_p)_\ell(\pa_{p+1}(\varphi_{p+1}.[B_{p+1}(M),B_{p}(M)])) &\neq& 0.
\eeq
We now introduce the following relative homological diagram.

\n Since we have supposed that $J$ has no critical points in $\Sig^+$, we derive from Proposition \ref{prop3.1} that the natural injection $\ov{i}\colon\;E_p\hookrightarrow F_{p+1}$ is an homology equivalence. Thus it induces a relative homology isomorphism, denoted
\beqq
(\ov{i})_\ell\colon\;H_\ell(E_{p+1},E_p) &\to& H_\ell(E_{p+1},F_{p+1}),\;\ell\geq 0.
\eeqq
Denote,
\beqq
\d_{p+1}\colon\;H_{\ell+1}(E_{p+1},E_p) &\to& H_\ell(E_{p},F_{p}),\;\ell\geq 0,
\eeqq
the connecting homomorphism of the triad $(E_{p+1},E_p,F_p)$. the following diagram commutes.
\[\xymatrix@C=3cm@R=1cm{
H_{\ell+1}(B_{p+1}(M),B_{p}(M)) \ar[r]^{(\phi_{p+1})_\ell} \ar[d]^{\pa_{p+1}}  & H_{\ell+1}(E_{p+1},F_{p+1}) \ar[d]^{\ov{\pa_{p+1}}} \\
H_{\ell}(B_{p}(M),B_{p-1}(M)) \ar[r]^{(\phi_{p})_\ell}  & H_{\ell+1}(E_{p+1},F_{p+1})
}\]
where $\ov{\pa_{p+1}}=\d_{p+1}o(\ov{i})_{\ell+1}^{-1}$. We then have
\[(\phi_{p})_\ell o\pa_{p+1}=\ov{\pa_{p+1}}o(\phi_{p+1})_\ell.\]
Thus from \eqref{eq3.15}, we have
\beqq
\ov{\pa_{p+1}}((\phi_{p+1})_\ell(\varphi_{p+1}.[B_{p+1}(M),B_p(M)])) &\neq& 0
\eeqq
and therefore,
\beqq
(\phi_{p+1})_\ell(\varphi_{p+1}.[B_{p+1}(M),B_p(M)]) &\neq& 0
\eeqq
Using the fact that $(\phi_{p+1})_\ell$ is $H^\ell(\O^{p+1}/_{\varrho})$-linear, we get
\beqq
\varphi_{p+1}.(\phi_{p+1})_\ell([B_{p+1}(M),B_p(M)]) &\neq& 0,
\eeqq
and hence
\beqq
(\phi_{p+1})_\ell([B_{p+1}(M),B_p(M)]) &\neq& 0.
\eeqq
Claim \eqref{eq3.14} is valid and the result follows.
\end{pf}

\n We now state the following result
\begin{proposition}\label{prop3.5}
There exists a positive integer $p_0$ large enough such that for any $p\geq p_0$,
\beqq
\phi_{p}(B_{p+1}(M),B_{p-1}(M)) &\subset& (F_p,F_p).
\eeqq
\end{proposition}
\begin{pf}
Let $u=\underset{i=1}{\overset{p}{\sum}}\a_i\d_{a_i}\in B_p(M)$.
\beqq
\phi_{p}(u) &=& \frac{\sum_{i=1}^{p}\ov{\a_i}PU_{(\ov{a_i},\l)}}{\abs{\sum_{i=1}^{p}\ov{\a_i}PU_{(\ov{a_i},\l)}}},
\eeqq
where $(\ov{a_1},\cdots,\ov{a_p},\ov{\a_1},\cdots,\ov{\a_p})=(h_p oj_p oi_po\pi_p)(u)$. Here $h_p$, $j_p$, $i_p$ and $\pi_p$ are defined in the proof of Proposition \ref{prop3.3}. By \eqref{eq3.4}, we have
\beqq
\phi_{p}(u) &\in& W_u(V(p,\e_p)).
\eeqq
It remains to prove that $\phi_{p}(u)\in J_{p^{\soexp{2}-1}\wtilde{S}_{HL}}$. For this we use the expansion of Proposition \ref{prop2.5} and the fact that
\beqq
H(a,a) &\leq& C,\;\;\forall a\in M,\\
G(a,b) &\geq& \g>0,\;\;\forall a,b\in M,\;\text{ such that  }\diff(a,b)\geq\diff_1,
\eeqq
we get
\beqq
J(\phi_p(u)) &\geq& P^{\soexp{2}-1}\wtilde{S_{HL}}(1+(\ov{c}-p\ov{\g}\frac{1}{\l^{n-2}}),
\eeqq
where $\ov{c}$ and $\ov{\g}$ are two position constants. Let $p_0$ such that $\ov{c}-p_0\ov{\g}<0$, we derive our result for any $p\geq p_0$.
\end{pf}

\n \textbf{Proof of Theorem \ref{th1.1} completed:} Let $\O\subset\Real^n$, $n\geq 3$, be a bounded domain such that $H_{k_0}(\O)\neq 0$, $k_0\geq 1$. there exists a singular $k_0$-chain in $\O$ without boundary with $\Z_2$-coefficients, ($k_0$-cycle), which defines a non zero class in the homology group $H_\ell(\O)$, $\ell\geq 0$. following \cite{4 BB}, the $k_0$-cycle can be viewed as a compact and closed $k_0$-dimensional manifold $M$ in $\O$. Using the preceding notations, we prove that
\beq\label{eq3.16}
(\phi_{1})_\ell([B_{1}(M),B_0(M)]) &\neq& 0.
\eeq
Indeed, if $(\phi_1)_\ell([B_1(M),B_0(M)])=0$, then using the mapping $\chi$ introduced in \eqref{eq3.3}, we find that
\beqq
(\chi)_\ell((\phi_1)_\ell([B_1(M),B_0(M)])) &=& [M,\pa M]=0\;\;\text{ in  }H_\ell(\O).
\eeqq
This is absurd, since $M$ has no boundary and defines a non zero class in $H_\ell(\O)$. The assertion \eqref{eq3.16} holds and hence claim \eqref{eq1.5} of the first step of the introduction is valid. Now if we suppose that $J$ has no critical point in $\Sig^+$, it is proved in \eqref{eq3.13} and  \eqref{eq3.14} that for any $p\in\Natl\backslash\{0\}$,
\beqq
(\phi_{p})_\ell([B_{p}(M),B_{p-1}(M)])\neq 0 &\Rightarrow& (\phi_{p+1})_\ell([B_{p+1}(M),B_{p}(M)])\neq 0.
\eeqq
Thus \eqref{eq1.6} follows. Lastly, it is proved in Proposition \ref{prop3.5} that there exists a large positive integer $p_0$ such that the mapping $\phi_{p_0}$ is valued in the pair $(F_{p_0},F_{p_0})$. Therefore it is homologically trivial. This yields the assertion \eqref{eq1.7} and proves our Theorem.


\begin{thebibliography}{99}
\bibitem{BC 1991} A. Bahri and J. M. Coron, The scalar curvature problem on the standard three dimensional spheres, J. Funct. Anal. 95 (1991), 106–172.
\bibitem{BC 1989} A. Bahri and J. M. Coron, On a nonlinear Elliptic equation Involving the critical Sobolev Exponent: The effect of the topology on the domain, Comm. Pure Appl. Math. 41 (1988), 253–294.
\bibitem{BL} A. Bahri, Critical Point at Infinity in Some Variational Problems, Pitman Res. Notes Math., vol. 182, Longman Sci. Tech., Harlow, 1989. 81.
\bibitem{P2}E. de S. Böer, O.H. Miyagaki and P. Pucci, Existence and multiplicity results for a class of Kirchhoff-Choquard equations with a generalized sign-changing potential, Atti Accad. Naz. Lincei Rend. Lincei Mat. Appl., Special Issue: dedicated to the memory of Professor Antonio Ambrosetti, 33 (2022), 651-675.
  
    
\bibitem{2 BB} G. Bredon, Introduction to compact transformations Groups, Academic Press, New York (1972).

\bibitem{BN} H. Brézis and L. Nirenberg, Positive solutions of nonlinear elliptic equations involving critical Sobolev exponents, Comm. Pure Appl. Math., 36 (1983), 437–477.



\bibitem{16 Aut} D. Cassani, J. V. Scahftingen and J. Zhang, Groundstates for Choquard type equations with Hardy-Littlewood-Sobolev lower critical exponent, Available via https://arxiv.org/pdf/1709.09448.pdf.
    
    \bibitem{C} J. M. Coron, Topologie et cas limite des injections de Sobolev, C. R. Acad. Sci. Paris Sér. I Math. 299 (1984), no. 7, 209–212.

\bibitem{4 BB} A. Dold, Lectures on Algebraic Topology, Springers-Verlag Berlin and New York.

\bibitem{17 Pd} L. Du and M. Yang, Uniqueness and nondegeneracy of solutions for a critical nonlocal equation, Discrete Contin. Dyn. Syst. 39, 5847–5866 (2019)
    
    \bibitem{19 Green} F. Gao and M. Yang, On nonlocal Choquard equations with Hardy–Littlewood–Sobolev critical exponents, J. Math. Anal. Appl., 448 (2017), 1006–1041.

\bibitem{20 Pd} F. Gao and M. Yang, The Brezis–Nirenberg type critical problem for the nonlinear Choquard equation, Sci. China Math. 61, 1219–1242 (2018).
\bibitem{31 Aut} D. Goel, V. Radulescu and K. Sreenadh, Coron problem for nonlocal equations involving Choquard nonlinearity, Adv. Nonlinear Stud. (2019), https://doi.org/10.1515/ans-2019-2064.

\bibitem{11 2023} D. Goel and K. Sreenadh, Critical growth elliptic problems involving Hardy-Littlewood-Sobolev critical exponent in non-contractible domains, Adv. Nonlinear Anal. 9 (2020), no. 1, 803–835.

\bibitem{23 Pd} L. Guo, T. Hu, S. Peng ans W. Shuai, Existence and uniqueness of solutions for Choquard equation involving Hard–Littlewood–Sobolev critical exponent, Calc. Var. Partial Differ. Equ. 58, 128, 34 pp (2019),

\bibitem{P1} S. Liang, P. Pucci and B. Zhang, Multiple solutions for critical Choquard-Kirchhoff type equations, Adv. Nonlinear Anal. 10 (2021), 400-419.
\bibitem{P3} S. Liang, P. Pucci and Y. Song, On the critical Choquard-Kirchhoff p-sub-Laplacian equation in the entire Heisenberg group, Anal. Geom. Metr. Spaces, Special Issue: Second Order Subelliptic PDEs dedicated to the 80th burthday of Professor Ermanno Lanconelli, 12 (2024),Paper No. 20240006, pages 28
\bibitem{15 DCDS} E. Lieb, Sharp constants in the Hardy-Littlewood-Sobolev and related inequalities, Ann. of Math. 118(1983), 349–374.
\bibitem{21 BN} E. Lieb, Existence and uniqueness of the minimizing solution of Choquard’s nonlinear equation, Studies in Appl. Math., 57(1976/77), 93–105.
\bibitem{17 DCDS} E. Lieb and M. Loss, "Analysis". Gradute Studies in Mathematics, AMS, Providence, Rhode island, 2001.
\bibitem{Lions} P. L.  Lions, The concentration compactness principle in the calculus of variations, (Part 1 and Part 2), Riv. Mat. Iberoamericana 1 (1985), 145-201, 45–121.
\bibitem{19 DCDS} P.L. Lions, The Choquard equation and related questions, Nonlinear Anal. 4 (1980), 1063–1072.

\bibitem{25 DCDS} L. Ma and L. Zhao, Classification of positive solitary solutions of the nonlinear Choquard equation, Ration. Mech. Anal. 195(2010), 455–467.

\bibitem{23 2023} V. Moroz and J. V. Schaftingen, A guide to the Choquard equation, J. Fixed Point Theory Appl. 19 (2019), 773–813.

\bibitem{47 Aut} V. Moroz and J. V. Schaftingen, Groundstates of nonlinear Choquard equations: HardyLittle-woodSobolev critical exponent, Commun. Contemp. Math. 17, (2015).
\bibitem{19 2023} V. Moroz and J. V. Schaftingen, Groundstates of nonlinear Choquard equations: existence, qualitative properties and decay asymptotics, J. Funct. Anal. 265 (2013), no. 2, 153–184,

\bibitem{29 BN} S. Pekar, Untersuchungüber die Elektronentheorie der Kristalle, Akademie Verlag, Berlin, 1954.
\bibitem{30 BN} R. Penrose, On gravity’s role in quantum state reduction, Gen. Relativ. Gravitat. 28 (1996), 581–600.
\bibitem{SYZ Pdelta} M. Squassina, M. Yang and S. Zhao, Local uniqueness of blow-up solutions for critical Hartree equations in bounded domain,  Calc. Var. (2023) 62:217, https://doi.org/10.1007/s00526-023-02551-1.
\bibitem{Struwe} M. Struwe, A global compactness result for elliptic boundary value problem involving limiting nonlinearities, Math. Z. 187, (1984), 511-517.

\end{thebibliography}
\end{document}